\documentclass[12pt,a4paper]{article}

\usepackage{amsmath}    
\usepackage{amssymb}
\usepackage{graphicx}   
\usepackage{verbatim}   
\usepackage{color}      
\usepackage{hyperref}   
\usepackage{enumerate}
\usepackage{listings}
\usepackage{color}
\usepackage{slashed}
\usepackage{amsthm}
\usepackage{bbold}
\usepackage{listing}
\usepackage{mathrsfs}
\usepackage{empheq}
\usepackage[british]{babel}
\usepackage[font={small,it}]{caption}
\usepackage{etoolbox}

\usepackage{cleveref}
\crefformat{footnote}{#2\footnotemark[#1]#3}

\definecolor{mygreen}{rgb}{0,0.6,0}
\definecolor{mygray}{rgb}{0.5,0.5,0.5}
\definecolor{mymauve}{rgb}{0.58,0,0.82}

\AtEndEnvironment{example}{%
  \qed
}

\newcommand{\del}{\partial}
\renewcommand{\phi}{\varphi}
\newcommand{\vecc}[2]{\left ( \begin{array}{c}#1\\#2\\ \end{array}\right )}
\newcommand{\veccc}[3]{\left ( \begin{array}{c}#1\\#2\\#3\\ \end{array}\right )}
\newcommand{\dd}{\mathrm{d}}
\newcommand{\ee}{\text{e}}
\newcommand{\ii}{\mathbb{i}}
\newcommand{\id}{\mathbb{1}}

\renewcommand{\and}{\wedge}

\renewcommand{\vec}{\mathbf}

\newtheorem{theorem}{Theorem}[section]
\newtheorem{corollary}{Corollary}[section]

\newtheorem{definition}{Definition}[section]

\newtheorem{remark}{Remark}[section]

\renewcommand{\title}{An Active Flux method for the Euler equations based on the exact acoustic evolution operator}

\newcommand{\authorOne}{Wasilij Barsukow\footnote{CNRS, Institut de Mathématiques de Bordeaux (IMB), UMR 5251, 351 Cours de la Libération, 33405 Talence, France, wasilij.barsukow@math.u-bordeaux.fr}}
\begin{document}

\begin{center} \Large
\title

\vspace{1cm}

\date{}
\normalsize

\authorOne
\end{center}

\begin{abstract}

A new Active Flux method for the multi-dimensional Euler equations is based on an additive operator splitting into acoustics and advection. The acoustic operator is solved in a locally linearized manner by using the exact evolution operator. The nonlinear advection operator is solved at third order accuracy using a new approximate evolution operator. To simplify the splitting, the new method uses primitive variables for the point values and for the reconstruction. In order to handle discontinuous solutions, a blended bound preserving limiting is used, that combines a priori and a posteriori approaches. The resulting method is able to resolve multi-dimensional Riemann problems as well as low Mach number flow, and has a large domain of stability.

Keywords: Active Flux, Euler equations

Mathematics Subject Classification (2010): 65M08, 65M70, 76M12, 35L45

\end{abstract}

\section{Introduction}

The Active Flux method, originally introduced for one-dimensional (1-d) linear advection in \cite{vanleer77}, has received much attention since the pioneering works \cite{eymann11,eymann13}. As degrees of freedom, it employs averages and point values at cell interfaces. These are independent, i.e. the Active Flux method consists of update equations for both the averages and the point values. The former is easy to obtain for conservation laws: since the point values are located at cell interfaces, the flux can be directly evaluated; no Riemann solvers are needed and one recognizes the continuous nature of the spatial approximation. While the update of the averages is ``exact'', i.e. as accurate as the point values, it cannot incorporate upwinding, since the latter always implies some kind of additional, artificial diffusion.

It is the point values which need to include the upwinding necessary for stability. The initial method from \cite{vanleer77} uses characteristic tracing and a reconstruction, whose value at the foot of the characteristic is taken as the point value at the next time step. This procedure naturally includes upwinding, of course, and yields a one-stage method stable up to $\text{CFL} = 1$. Generalizations of this update procedure appeared under the name of \emph{(approximate) evolution operators} in e.g. \cite{eymann13,fan17,barsukow18activeflux,barsukow19activeflux,chudzik24}. 

It was shown in \cite{barsukow19activeflux} that one needs to go one order of accuracy beyond local linearization to obtain a third-order Active Flux method. For Burgers' equation $\del_t u + u \del_x u = 0$, instead of the local linearization $x \mapsto x - u(x) t$, for example, \cite{roe17} suggests to use
\begin{align}
 x \mapsto x - \frac{u(x) t}{1 + t \del_x u(x)} \simeq x - u(x) t + t^2 u(x) \del_x u(x) + \mathcal O(t^3)
\end{align}
while \cite{barsukow19activeflux} proposes to iterate
\begin{align}
 x \mapsto x - u\Big(x - u(x) t\Big) t \simeq x - u(x) t + t^2  u(x) \del_x u(x ) + \mathcal O(t^3)
\end{align}

In \cite{barsukow19activeflux}, a sufficiently accurate approximate evolution operator was achieved for any hyperbolic system of conservation laws in one spatial dimension. It has been shown that it is by no means sufficient to merely iterate the linearization a few times to achieve third order of accuracy -- this only works for scalar conservation laws since their characteristics are straight even in the nonlinear case (see Section \ref{sec:highorder} for further explanations). An alternative, ADER-inspired approach is \cite{kerkmann18}.

Concerning multiple spatial dimensions, early effort (\cite{fan17,barsukow18activeflux}) focused on the equations of linear acoustics and the corresponding evolution operator. The exact solution for a general class of data (in particular those not differentiable) was obtained in \cite{barsukow17}, initially destined to study genuinely multi-dimensional (multi-d) Godunov methods. The exact solution was used to demonstrate that even the complete solution of the multi-dimensional Riemann problem does not render a Godunov method stationarity or vorticity preserving. In \cite{barsukow18activeflux} the exact acoustic evolution operator was used to construct an Active Flux method on two-dimensional Cartesian grids. Contrary to the Godunov method, the Active Flux method did turn out to be stationarity and vorticity preserving.

To construct an evolution operator for nonlinear systems of conservation laws in multi-d is significantly more difficult, since third order accuracy requires going beyond local linearization. At this point, in \cite{abgrall20,abgrall22}, a semi-discrete version of Active Flux appeared. There, the point value is evolved according to an ODE in time, with the space derivative replaced by a finite difference formula (a non-standard one, since it involves both point values and averages in some neighbourhood of the degree of freedom). The time discretization is provided by a standard Runge-Kutta method according to the method of lines. The advantage of this approach is its immediate applicability to nonlinear systems in multi-d (e.g. \cite{abgrall24prim,barsukow24afeuler,duan24}), with the structure preserving properties still in place (\cite{barsukow24affourier}). A slight disadvantage are the reduced low CFL numbers (0.2 and less), which \cite{roe21} connects to the usage of a multistage Runge-Kutta method. In fact, particular choices of less compact finite difference formulas do yield stability all the way up to $\text{CFL} = 1$ for RK3 (\cite{abgrall22}), but they are difficult to generalize in multi-d. However, undeniably, a Runge-Kutta method adds its own numerical diffusion.

The continuous approximation lets the semi-discrete Active Flux method appear different from, say, the Discontinuous Galerkin method, but the absence of the evolution operator also moves it away from van Leer's and Roe's original ideas. In the context of implicit methods for linear advection it was shown in \cite{barsukow23afim} that one-stage implicit Active Flux methods obtained using an evolution operator were significantly better than semi-discrete Active Flux methods integrated with standard (implicit) integrators, at the same (or less) cost. This allows to conjecture that it might be worth striving for an Active Flux method for the Euler equations in multi-d that is based on evolution operators.

To the author's knowledge, two such evolution-operator based Active Flux methods for multi-d Euler equations are available in the literature.
\cite{roe18} suggests to start out from a Lax-Wendroff/Cauchy-Kowalevskaya expansion. The terms $\mathcal O(t)$ involve locally linearized acoustics and advection separately, while terms $\mathcal O(t^2)$ involve deviations from a pure local linearization and interactions between these operators. The treatment of these extra terms remains difficult. The different operators of the Euler equations were treated in \cite{fan17,maeng17}, yet a clear prescription how the problematic terms are to be discretized seems still to be missing, as is an investigation what impact particular choices would have in practice.

The recent work \cite{chudzik24} employs second-order approximations to the evolution operators of linear advection and acoustics, but does not split them. This requires treating all the ways how a tilted characteristic cone (subject to ``wind'') can cut through the cells in the vicinity of a point value. Here, the high-order terms that would be required for a third-order method are not included. 

In view of such difficulties there seems to exist a certain interest in accepting second-order accuracy in time in certain ingredients of the numerical method for the sake of computational efficiency. The practical value of a time integration whose order of accuracy is not quite matching that of space discretization might depend on the application in question and on the willingness to spend effort on implementation and to accept a possible increase in computational effort. It is, for example common practice to use e.g. RK3 in conjunction with arbitrarily high-order methods (e.g. \cite{jiang96,vilar19,abgrall22}). Formally second-order methods that employ space reconstruction of degree up to 6 can bring significant improvement even without the corresponding increase in order of accuracy in time (\cite{leidi24}).

This paper proposes to use
\begin{enumerate}[--]
 \item a simple operator split between acoustics and advection, unsurpassed in computational efficiency, yet second-order; with
 \item a fully third-order approximate evolution operator for the (nonlinear) advection in multi-d, and
 \item a locally linearized, yet otherwise exact evolution operator for acoustics.
\end{enumerate}
The advantage of the exact evolution operator is the inclusion of all the multi-dimensional information. Moreover, since the linearized method reduces to \cite{barsukow18activeflux}, its von Neumann stability limit is the same (contrary to \cite{chudzik24}). The present paper gives a large amount of detail concerning the efficient implementation of the exact evolution operator. A bound-preserving limiting strategy combining one similar to that in \cite{duan24} with an a posteriori limiting for the point values is presented as well.

The result is a method for the Euler equations, possibly having some similarity with \cite{roe18}, yet also clear differences, e.g. in the choice of a non-conservative advective operator and the usage of primitive variables for the reconstruction and the point values (cell averages are, naturally, in conserved variables). The simultaneous usage of primitive and conservative formulations has been advertised in \cite{abgrall24prim} for the semi-discrete Active Flux, and already in \cite{barsukow19activeflux} primitive variables at point values were in use for the fully discrete Active Flux method. The difference of the present work to \cite{chudzik24} is the usage of the exact evolution operator, with more accuracy at essentially the same computational cost due to efficient implementation, and the algorithmic simplification that is brought about by the operator splitting. The purpose of this paper therefore is to present this simple one-stage method for the multi-d Euler equations. Beside that, the implementation details concerning the exact evolution operator for linear acoustics, at that time not included in \cite{barsukow18activeflux}, hopefully can now serve as a reference for future use.

The paper's organization is as follows: Section \ref{sec:evoop} introduces the efficient implementation of the evolution operator for linear acoustics, complemented by Appendix \ref{sec:sphericalmean} detailing the computation of spherical means.
 The method for the Euler equations and its properties are described in Section \ref{sec:euler}, followed by the numerical results of Section \ref{sec:numerical}. The Appendices \ref{app:evoop} and \ref{app:details} contain pseudocode implementations of the acoustic evolution operator.

In the following, computations often are performed in three space dimensions, and two-dimensional results are obtained by assuming that the functions involved do not depend on $z$. Vectors associated with spatial directions and thus naturally having the same number of components as the dimension of the space are set in boldface with, in particular, $\vec x = (x,y,z)^\text{T}$. The notation $\vec a \cdot \vec b$ for $\vec a, \vec b \in \mathbb R^n$ denotes the Euclidean scalar product of $\vec a$ and $\vec b$; $S^2 := \{ \vec x \in \mathbb R^3 : |\vec x| = 1 \} $ is the unit sphere in three space dimensions.

\section{Evolution operator for linear acoustics} \label{sec:evoop}

The exact solution of linear acoustics
\begin{align}
 \del_t \vec v + c \nabla p &= 0  & \vec v &\colon \mathbb R^+_0 \times \mathbb R^3 \to \mathbb R^3\\
 \del_t p + c \nabla \cdot \vec v &= 0 & p &\colon \mathbb R^+_0 \times \mathbb R^3 \to \mathbb R
\end{align}
on $\mathbb R^3$ endowed with initial data $(\vec v_0, p_0)$ is
\begin{align}
 p(t, \vec x) &= \del_r \left(r M[p_0](\vec x, r) \right) - \frac{1}{r} \del_r \left(r^2 M[\vec n \cdot \vec v_0](\vec x, r)\right ) \label{eq:preformulfunc}\\
 \vec v(t, \vec x) &= \frac23 \vec{ v}_0(\vec x) - \frac{1}{r} \del_r \left( r^2 M[p_0  \vec n](\vec x, r) \right) + \del_r\left(rM[(\vec v_0 \cdot \vec n)  \vec n](\vec x, r) \right ) \nonumber \\&- M\left[  \vec v_0 - 3  (\vec v_0 \cdot \vec n)  \vec n   \right ](\vec x, r) - \int_0^{ct} \dd r \frac{1}{r} M \left[  \vec v_0 - 3  (\vec v_0 \cdot \vec n)  \vec n   \right ](\vec x, r) \label{eq:vsolution23func}
\end{align}
where all formulae are understood to be evaluated at $r=ct$ after performing all the differentiations (see \cite{barsukow17} for the derivation) and $M[f] \colon \mathbb R^3 \times \mathbb R^+_0 \to \mathbb R$ is the spherical mean of $f \colon \mathbb R^3 \to \mathbb R$, defined as
\begin{align}
 M[f](\vec x,r) &:= \frac{1}{4\pi} \int_{S^2} \dd \vec y \, f(\vec x + r \vec y) 
\end{align}
For the definition of spherical means such as $M[\vec n \cdot \vec v_0]$, see Section \ref{sec:sphericalmean}. The data $(\vec v_0, p_0)$ needs generally to be twice continuously differentiable for the solution formula from \cite{fan17} to be applicable, while \eqref{eq:preformulfunc}--\eqref{eq:vsolution23func} has the advantage that it only requires regularity with respect to the distance variable $r$. Thus, it remains valid even if the data are discontinuous, as long as the discontinuities are straight lines through $\vec x$, and this is precisely what is needed in the numerical setting. The formulation \eqref{eq:preformulfunc}--\eqref{eq:vsolution23func} thus avoids having to deal with distributions.

This formula is valid in three spatial dimensions, but one can of course restrict oneself to initial data that do not depend on certain variables, or initial velocity fields with certain components zero.
In the following, $u$ and $v$ will be used to denote the $x$ and $y$ components of $\vec v$.

The degrees of freedom of Active Flux are located at edges and nodes, with data being piecewise polynomial. It thus makes sense to divide the solution up in contributions of half-planes/quadrants. Consider, more generally, initial data in some wedge $W = [\phi_\text{min}, \phi_\text{max}) \times [0, \pi]$ centered at $\vec x$ and denote by $p^{W}(t, \vec x), \vec v^W(t, \vec x)$ the solution obtained by taking in \eqref{eq:preformulfunc}--\eqref{eq:vsolution23func} everywhere spherical averages $M^W$ over the wedge instead of $M$. Then, if $\{ W^i \}_i$ is a finite partition of $[0, 2 \pi) \times [0, \pi]$, then 
\begin{align}
 M[f](x, r) &= \frac1{4\pi} \int_{\bigcup_i W_i} \dd \vec y f(\vec x + r \vec y) = \frac1{4\pi} \sum_i 2 (\Delta \phi)_i \frac{1}{2(\Delta \phi)_i} \int_{W_i} \dd \vec y f(\vec x + r \vec y) \\
 &= \frac1{2\pi} \sum_i (\Delta \phi)_i  M^{W_i}[f](x, r)
\end{align}
Thus,
\begin{align}
 p(t, \vec x) &= \frac1{2\pi} \sum_i |W_i| p^{W_i}(t, \vec x) & \vec v(t, \vec x) &= \frac1{2\pi} \sum_i |W_i| \vec v^{W_i}(t, \vec x)
\end{align}

The ability to combine solutions across wedges as a simple sum is a consequence of the solution formula above. Finally, when applying the solution formula for every degree of freedom, $\vec x$ is fixed and the result is simply a polynomial in $t$ with easily precomputable coefficients. For the sake of giving the reader a quicker access to the Section \ref{sec:euler}, which goes on to develop a method for the Euler equations, detailed computations are found in the Appendix \ref{sec:sphericalmean}.

\section{An Active Flux method for the Euler equations} \label{sec:euler}

\subsection{General remarks}

\begin{figure}[h]
 \centering
 \includegraphics[width=0.25\textwidth]{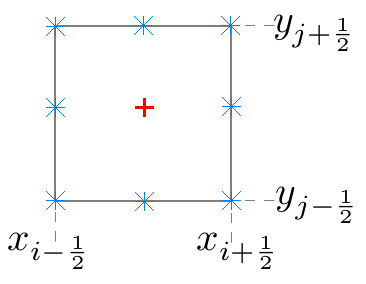}
 \caption{Degrees of freedom of Active Flux: Averages and point value at nodes and edge midpoints.}
 \label{fig:dof}
\end{figure}

On a Cartesian grid with cells
\begin{align}
 C_{ij} = \left[x_{i-\frac12},x_{i+\frac12} \right ] \times \left[ y_{j-\frac12}, y_{j+\frac12}  \right]
\end{align}
the degrees of freedom of Active Flux are (see also Figure \ref{fig:dof}) cell averages
\begin{align}
 \bar q_{ij}(t) \simeq \frac{1}{|C_{ij}|} \int_{C_{ij}} \dd \vec x \, q(t, \vec x)
\end{align}
and point values at nodes
\begin{align}
 q_{i+\frac12,j+\frac12}(t) \simeq q(t, x_{i+\frac12}, y_{j+\frac12})
\end{align}
and at edge midpoints
\begin{align}
 q_{i+\frac12,j}(t) &\simeq q(t, x_{i+\frac12}, y_{j}) & q_{i,j+\frac12}(t) &\simeq q(t, x_{i}, y_{j+\frac12})
\end{align}

Time is discretized at instants $t^n$, $n \in \mathbb N_0$, with $\Delta t$ a generic name for the difference $t^{n+1} - t^n$. The values of the degrees of freedom at instant $t^n$ are denoted by a superscript $n$.

A reconstruction in each cell serves as initial data for the evolution operator. The classical reconstruction in Active Flux agrees pointwise with the given point values, and its average equals the given cell average. With 9 conditions, a biparabolic reconstruction in
\begin{align}
\mathrm{span}(1, x, x^2, y, xy, x^2y, y^2, xy^2, x^2 y^2) 
\end{align}
is natural (as used in \cite{barsukow18activeflux,kerkmann18}, where the polynomial is given explicitly). Since the restriction of a biparabolic function on one of the (axis-aligned) edges is a parabola, and since it is uniquely determined by the three values present along each edge, the reconstructions form globally continuous initial data. As is described later, a biparabolic reconstruction is still used here, but computed slightly differently.

\subsection{Update of the average}

Recall that the update of the average $\bar q_{ij} := \frac{1}{|C_{ij}|} \int \dd \vec x \, q(t, \vec x)$ is obtained immediately by integrating the conservation law 
\begin{align}
 \del_t q + \nabla \cdot \vec f(q) &= 0 & q &\colon \mathbb R^+_0 \times \mathbb R^2 \to \mathbb R^m \\
 && \vec f &= (f^x, f^y) \quad f^x, f^y \colon \mathbb R^m \to \mathbb R^m
\end{align}
over the cell $C_{ij}$ and over the time step $[t^n, t^{n+1}]$ and replacing the integrals by space-time Simpson-rule quadratures (see e.g. \cite{eymann13}): 
\begin{align}
 \frac{\bar q_{ij}^{n+1} -\bar q_{ij}^n}{\Delta t} &=-\frac{1}{\Delta x}\left( \hat f^x_{i+\frac12,j} - \hat f^x_{i-\frac12,j} \right) -\frac{1}{\Delta y}\left( \hat f^y_{i,j+\frac12} - \hat f^y_{i,j-\frac12} \right)  \simeq - \int_{\del C_{ij}} \dd \vec x \, \vec n \cdot \vec f(q)
 \end{align}
with {\footnotesize 
\begin{align}
 \hat f^x_{i+\frac12,j} &:= \frac16 \left( \frac{f^x(q^n_{i+\frac12,j-\frac12})}{6} + \frac{2f^x(q^n_{i+\frac12,j})}{3} + \frac{f^x(q^n_{i+\frac12,j+\frac12})}{6}  \right) \label{eq:fluxquadraturex}\\ 
 &\!\!\!\!\!\!\!\!\!\!\!\!\!\!\!\!+ \nonumber  \frac23 \left( \frac{f^x(q^{n+\frac12}_{i+\frac12,j-\frac12})}{6} + \frac{2f^x(q^{n+\frac12}_{i+\frac12,j})}{3} + \frac{f^x(q^{n+\frac12}_{i+\frac12,j+\frac12})}{6}  \right)
 +\frac16 \left( \frac{f^x(q^{n+1}_{i+\frac12,j-\frac12})}{6} + \frac{2f^x(q^{n+1}_{i+\frac12,j})}{3} + \frac{f^x(q^{n+1}_{i+\frac12,j+\frac12})}{6}  \right) \\
 &\simeq \frac{1}{\Delta y} \int_{y_{j-\frac12}}^{y_{j+\frac12}} \dd y \frac{1}{\Delta t} \int_{t^n}^{t^{n+1}} \dd t \, f^x(q(x_{i+\frac12}, y, t)) \\
 \hat f^y_{i,j+\frac12} &:= \frac16 \left( \frac{f^y(q^n_{i-\frac12,j+\frac12})}{6} + \frac{2f^y(q^n_{i,j+\frac12})}{3} + \frac{f^y(q^n_{i+\frac12,j+\frac12})}{6}  \right) \\ 
 &\!\!\!\!\!\!\!\!\!\!\!\!\!\!\!\!+ \nonumber  \frac23 \left( \frac{f^y(q^{n+\frac12}_{i-\frac12,j+\frac12})}{6} + \frac{2f^y(q^{n+\frac12}_{i,j+\frac12})}{3} + \frac{f^y(q^{n+\frac12}_{i+\frac12,j+\frac12})}{6}  \right)
 +\frac16 \left( \frac{f^y(q^{n+1}_{i-\frac12,j+\frac12})}{6} + \frac{2f^y(q^{n+1}_{i,j+\frac12})}{3} + \frac{f^y(q^{n+1}_{i+\frac12,j+\frac12})}{6}  \right) \\ &\simeq \frac{1}{\Delta x} \int_{x_{i-\frac12}}^{x_{i+\frac12}} \dd x \frac{1}{\Delta t} \int_{t^n}^{t^{n+1}} \dd t \, f^y(q(x, y_{j+\frac12}, t))
\end{align}}

Since this step requires the new point values at $t^{n} + \frac{\Delta t}{2}$ and $t^{n+1} = t^n + \Delta t$, it is performed as the very last. 

\subsection{Update of the point values}

Recall the Euler equations in primitive variables
\begin{align}
 \del_t \rho + \vec v \cdot \nabla \rho + \rho \nabla \cdot \vec v &= 0 \\
 \del_t \vec v + (\vec v \cdot \nabla ) \vec v + \frac{\nabla p}{\rho} &= 0 \\
 \del_t p + \vec v \cdot \nabla p + \rho c^2 \nabla \cdot \vec v &= 0
\end{align}

Contrary to e.g. \cite{barsukow24afeuler} here the following special choices are made:
\begin{enumerate}
\item The first ingredient of the presented approach is the usage of primitive variables for the point values. 
\item The second is a reconstruction in primitive variables. 
For point values, the transformation between conserved and primitive variables is immediate. Therefore, here, after transforming the point values to conserved variables, and given the average over the cell (in conserved variables, too), a two-dimensional Simpson rule is inverted to find the value at cell center (in conserved variables). This latter is finally transformed to primitive variables and a biparabolic reconstruction is computed. 
\item An operator splitting is employed, by first solving
\begin{align}
 \del_t \rho + \rho \nabla \cdot \vec v &= 0 \\
 \del_t \vec v + \frac{\nabla p}{\rho} &= 0 \\
 \del_t p + \rho c^2 \nabla \cdot \vec v &= 0
\end{align}
in a locally linearized manner, i.e. actually solving
\begin{align}
 \del_t \vec v + c_0 \nabla \left(\frac{p}{\rho_0c_0}\right) &= 0 \label{eq:linacousticsfromeuler1}\\
 \del_t \left(\frac{p}{\rho_0c_0}\right) + c_0 \nabla \cdot \vec v &= 0\label{eq:linacousticsfromeuler2}
\end{align}
with $c_0$, $\rho_0$ evaluations of the speed of sound and of the density at the location of the respective point value, and obtaining the density from
\begin{align}
\del_t (\rho c_0^2 - p) &= 0 
\end{align}
To solve \eqref{eq:linacousticsfromeuler1}--\eqref{eq:linacousticsfromeuler2}, the exact evolution operator for linear acoustics is used, as described in Section \ref{sec:evoop}. The biparabolic reconstruction (involving values at time $t^n$) serves as the initial datum.

In a second step,
\begin{align}
 \del_t \rho + \vec v \cdot \nabla \rho  &= 0 \label{eq:advectionfromeuler1}\\
 \del_t \vec v + (\vec v \cdot \nabla ) \vec v  &= 0 \label{eq:advectionfromeuler2}\\
 \del_t p + \vec v \cdot \nabla p  &= 0\label{eq:advectionfromeuler3}
\end{align}
is solved. Local linearization would amount to evolving the data according to 
\begin{align}
\vec x \mapsto \vec x - \vec v_0 t \label{eq:characteristicsimple}
\end{align}
with, again, $\vec v_0$ the evaluation of the velocity at the location of the respective point value. Here instead, a fully third-order algorithm is used to solve \eqref{eq:advectionfromeuler1}--\eqref{eq:advectionfromeuler3}, which is discussed in Section \ref{sec:highorder}. The same biparabolic reconstruction (involving values at time $t^n$) serves as the initial datum. 

\item The two operators are combined in \emph{additive} operator splitting.

\end{enumerate}

\begin{remark}
 Define solution operators $S_1, S_2 \colon \mathbb R^M \times [0, T] \to \mathbb R^M$, such that $S_k(q_0, t)$, $k=1,2$ is the solution at time $t$ of $\text{IVP}_k$ (obtained by splitting some IVP into two operators) with initial data $q_0 \in \mathbb R^M$. Then, recall that multiplicative operator splitting is 
 \begin{align}
  S_2(S_1(q_0, t), t)
 \end{align}
 while 
 \begin{align}
  S_2(q_0, t) + S_1(q_0, t) - q_0
 \end{align}
 is additive operator splitting. Multiplicative operator splitting would require the solution of the first sub-operator to serve as initial data of the second. However, it is not possible to update the average individually with any of the two sub-operators that are considered here, since the presented splitting into acoustics and advection yields non-conservative sub-operators. Additive operator splitting can be applied to the same initial data and allows to circumvent this difficulty. The update of the average then happens using the unsplit equations.
\end{remark}

\begin{remark}
The cell average is evolved in conserved variables by evaluating the flux $$\vec f(q_\text{prim}) = \veccc{\rho \vec v}{\rho \vec v \otimes \vec v + p \id}{\vec v \left(\frac{\gamma p}{\gamma-1} + \frac12 \rho |\vec v|^2\right)}$$ using the primitive variables $(\rho, \vec v, p)$ stored at the point values. The method thus is conservative.
\end{remark}

\begin{remark}
 The transformation of the point values back and forth between primitive and conserved variables does not entail a big cost. Using primitive variables for the reconstruction is a common approach for Finite Volume methods. More importantly, this extra effort is accepted here for the greater advantage of a significantly simpler usage of the acoustic evolution operator. In principle, however, it seems possible to perform essentially the same algorithm in conserved variables alone.
\end{remark}

Due to the usage of the exact evolution operator for acoustics, the linearization of the present scheme amounts to a combination of advection and acoustics, the stability of which has been analyzed in \cite{chudzik21}. The shortest separation between two degrees of freedom is $\mathrm{min}(\Delta x, \Delta y)/2$. The method thus has a maximal CFL number of 0.5.

\subsection{Order of accuracy} \label{sec:highorder}

\subsubsection{Overview}

It has been shown in \cite{barsukow19activeflux} that local linearization (i.e. the evaluation of the Jacobians of the sub-operators \eqref{eq:linacousticsfromeuler1}--\eqref{eq:linacousticsfromeuler2}, \eqref{eq:advectionfromeuler1}--\eqref{eq:advectionfromeuler3} at the point value that one seeks to update) leads to a second-order accurate method only. In the same work it has been shown that fixpoint iteration on the characteristic equation leads to high-order results for scalar conservation laws. This is because even if they are nonlinear, characteristics remain straight lines. This is no longer true for systems, and \cite{barsukow19activeflux} showed that a naive fixpoint iteration does \emph{not} improve the order beyond second. 

For systems, the characteristics are curved since their speed depends on all characteristic variables\footnote{Riemann invariants are sometimes mistaken for characteristic variables. Riemann invariants have the same value across a simple wave (e.g. contact discontinuity or rarefaction) and they exist for e.g. the Euler equations. Characteristic variables $Q$ are those for which the hyperbolic system can be written in the form $\del_t Q_k + \lambda_k(Q_1, \ldots, Q_N)\del_x Q_k = 0$. If they exist, they remain constant on ``their'' characteristic, which is a stronger statement not restricted to simple waves, but holding true everywhere in the region of smoothness of the solution. However, they exist only if an integrability condition is verified (the matrix $\del Q/\del q$ needs to be commuted with the derivative $\del_x$). For linear systems, this integrability condition is void, and hence generically across a $k$-wave all characteristic variables but the $k$-th are Riemann invariants. This is also the case for nonlinear problems if the integrability condition is true (any $2\times2$ system, for example, which includes shallow water and isentropic Euler). For the full Euler equations, the integrability condition is false, and characteristic variables do not exist. The explanation therefore is to be understood in the simplified setting of hyperbolic systems with characteristic variables.}, but only one of them is constant along the characteristic. The speed is thus changing in function of the values brought in by the other characteristic families. 

The correct way to achieve a third-order approximate evolution operator in one spatial dimension is thus characterized by a mixing of eigenvalues. For example, for a $2 \times 2$ system 
\begin{align}
 \del_t q_1 + \lambda_1(q_1 , q_2) \del_x q_1 &= 0 \label{eq:simplehyp1}\\
 \del_t q_2 + \lambda_2(q_1 , q_2) \del_x q_2 &= 0 \label{eq:simplehyp2}
\end{align}
with initial data $q_k(0, x) = q_{0,k}(x)$, local linearization produces $x -t \lambda_1\big(q_{1,0}(x), q_{2,0}(x)\big)$ as the estimate of the foot of the 1-characteristic passing through $x$ at time $t$. The estimate to next order of accuracy is (\cite{barsukow19activeflux})
\begin{align}
 x - t \lambda_1 \left( q_{1,0}\big( x - t \lambda_1(x)  \big), q_{2,0}\left( x - t \frac{\lambda_1(x) + \lambda_2(x)}{2}    \right)      \right)
\end{align}
Here, $\lambda_k(x)$ is short-hand for $\lambda_k(q_{1,0}(x), q_{2,0}(x))$. Simple characteristic-wise fixpoint iteration will never produce signals traveling at the average speed $\frac{\lambda_1(x) + \lambda_2(x)}{2}$.
For more complicated hyperbolic systems than \eqref{eq:simplehyp1}-\eqref{eq:simplehyp2} one additionally needs to take into account the changing eigenvectors, see \cite{barsukow19activeflux} for the general algorithm in 1-d. Therein, this algorithm has been successfully applied to the Euler equations.

For nonlinear systems in multiple spatial dimensions the question of high-order evolution operators beyond local linearization is open. In practice, local linearization is employed (e.g. in \cite{fan15,chudzik24}), thus formally reducing the order to second. The nonlinear version of advection, however, is amenable to higher-order evolution as discussed in \cite{maeng17}. Here, a different and inexpensive modification for solve \eqref{eq:advectionfromeuler1}--\eqref{eq:advectionfromeuler2} is used, as described below. 

Overall, since the reconstruction and all quadratures are third-order, the presented method is third-order accurate in space. To yield an Active Flux method that is also third-order accurate in time, the point updates have to have an error of at most $\mathcal O(t^3)$ (\cite{barsukow19activeflux}). While the advective operator is solved with this accuracy (see below), the locally linearized acoustics operator and the operator splitting both have errors $\mathcal O(t^2)$ and thus make the method formally second-order in time. This, at one order higher, finds parallels in the current practice: For high-order methods such as WENO or DG, it is actually customary to save on the effort of an adequately high-order time integration and a method like RK3 is often used. As is shown below, experimental results on the presented method give rise to the conclusion that the global error is not dominated by its second-order component in typical numerical tests. In fact, basing one's conclusion on certain experimental results one could easily mistake the method for a third-order accurate one.

\subsubsection{A high-order evolution operator for nonlinear advection}

Consider a hyperbolic system 
\begin{align}
\del_t q + J^x(q) \del_x q + J^y (q) \del_y q = 0 \label{eq:generalsystem}
\end{align}
with initial data $q(0, \vec x) = q_0(\vec x)$. We denote the components of the initial data as $q_{i,0}(\vec x)$.

\begin{theorem} \label{thm:advectionapproxevoop}
 If $J^x = v^x \id$, $J^y = v^y \id$ and both $v^x$ and $v^y$ are components of $q$ (as in \eqref{eq:advectionfromeuler1}--\eqref{eq:advectionfromeuler2}), then the solution of \eqref{eq:generalsystem} (if smooth) is
\begin{align}
 q(t, \vec x) = q_\text{loc.lin.}(t, \vec x) +  t^2 \nabla q_{0} \cdot ((\vec v_0 \cdot \nabla) \vec v_0)  + \mathcal O(t^3)
\end{align}
where $q_\text{loc.lin.}(t, \vec x)$ is obtained by considering $J^x$, $J^y$ fixed with whatever value they have at $\vec x$ at $t=0$, and $\vec v_0 = (v^x(0,\vec x), v^y(0,\vec x))$.
\end{theorem}
\begin{proof}
 We find (summation over repeated indices):
 \begin{align}
  \del_t q_i(t, \vec x) &= -J^x_{ij} \del_x q_j - J^y_{ij} \del_y q_j \\
  \del_t^2 q_i(t, \vec x) &= -J^x_{ij} \del_x \del_t q_j - J^y_{ij} \del_y \del_t q_j -\del_t J^x_{ij} \del_x q_j - \del_t J^y_{ij} \del_y q_j \\
  &=J^x_{ij} \del_x (J^x_{jk} \del_x q_k + J^y_{jk} \del_y q_k ) + J^y_{ij} \del_y (J^x_{jk} \del_x q_k + J^y_{jk} \del_y q_k ) \\
  &\nonumber \qquad + \frac{\del J^x_{ij}}{\del q_k} (J^x_{k\ell} \del_x q_\ell + J^y_{k\ell} \del_y q_\ell ) \del_x q_j + \frac{\del J^y_{ij}}{\del q_k} (J^x_{k\ell} \del_x q_\ell + J^y_{k\ell} \del_y q_\ell ) \del_y q_j \\
  &=J^x_{ij}  \left(J^x_{jk} \del_x^2 q_k + J^y_{jk} \del_x\del_y q_k  + \frac{ \del J^x_{jk}}{\del q_\ell} \del_x q_\ell \del_x q_k + \frac{\del J^y_{jk}}{\del q_\ell} \del_x q_\ell \del_y q_k \right) \\
  &\nonumber \qquad + J^y_{ij} \left(J^x_{jk} \del_x \del_y q_k + J^y_{jk} \del_y^2 q_k + \frac{\del J^x_{jk}}{\del q_\ell} \del_yq_\ell \del_x  q_k + \frac{\del J^y_{jk}}{\del q_\ell} \del_y q_\ell \del_y q_k \right) \\
  &\nonumber \qquad + \frac{\del J^x_{ij}}{\del q_k} (J^x_{k\ell} \del_x q_\ell + J^y_{k\ell} \del_y q_\ell ) \del_x q_j + \frac{\del J^y_{ij}}{\del q_k} (J^x_{k\ell} \del_x q_\ell + J^y_{k\ell} \del_y q_\ell ) \del_y q_j
 \end{align}
 The 4 terms not involving the derivatives of $J^x$, $J^y$ are those appearing for a linear problem and hence also those produced by local linearization. The remaining terms are (renaming $j \leftrightarrow k$ in the last parts)
 \begin{align}
  S_i &:= J^x_{ij}  \left( \frac{ \del J^x_{jk}}{\del q_\ell} \del_x q_\ell \del_x q_k + \frac{\del J^y_{jk}}{\del q_\ell} \del_x q_\ell \del_y q_k \right) 
  + J^y_{ij} \left( \frac{\del J^x_{jk}}{\del q_\ell} \del_yq_\ell \del_x  q_k + \frac{\del J^y_{jk}}{\del q_\ell} \del_y q_\ell \del_y q_k \right) \\
  &\nonumber \qquad + \frac{\del J^x_{ik}}{\del q_j} (J^x_{j\ell} \del_x q_\ell + J^y_{j\ell} \del_y q_\ell ) \del_x q_k + \frac{\del J^y_{ik}}{\del q_j} (J^x_{j\ell} \del_x q_\ell + J^y_{j\ell} \del_y q_\ell ) \del_y q_k \\
  &= \left( J^x_{ij}  \frac{ \del J^x_{jk}}{\del q_\ell}  
  + \frac{\del J^x_{ik}}{\del q_j} J^x_{j\ell} \right ) \del_x q_\ell \del_x q_k 
  + \left( J^x_{ij}  \frac{\del J^y_{jk}}{\del q_\ell}  
  + J^y_{ij} \frac{\del J^x_{j\ell}}{\del q_k}  
  + \frac{\del J^x_{i\ell}}{\del q_j} J^y_{jk}  
  + \frac{\del J^y_{ik}}{\del q_j} J^x_{j\ell} \right) \del_x q_\ell \del_y q_k 
  \\
  &\nonumber \qquad
  + \left( J^y_{ij} \frac{\del J^y_{jk}}{\del q_\ell}   
  + \frac{\del J^y_{ik}}{\del q_j} J^y_{j\ell} \right) \del_y q_\ell \del_y q_k
 \end{align}
 Consider now the special case $(J^x)_{ij} = v^x \delta_{ij}$, $(J^y)_{ij} = v^y \delta_{ij}$:
 \begin{align}
  S_i &= 2 v^x   \frac{  \del v^x }{\del q_\ell}  \del_x q_\ell \del_x q_i
  + \left( 
    2 v^x  \frac{\del v^y}{\del q_\ell}  \del_x q_\ell \del_y q_i 
  + 2 v^y \frac{\del v^x}{\del q_k}  \del_x q_i \del_y q_k \right)
  +  2 v^y \frac{\del v^y}{\del q_\ell}   \del_y q_\ell \del_y q_i \\
  &= 2 (v^x   \del_x v^x  +  v^y   \del_y v^x )\del_x q_i 
  + 2 (v^x    \del_x v^y +  v^y  \del_y v^y) \del_y q_i \\
  &= 2 \nabla q_i \cdot ((\vec v \cdot \nabla) \vec v)
 \end{align}
At this point the advantage of primitive variables is clearly visible. The statement is obtained by realizing that the derivatives need to be evaluated at $t=0$ in the Taylor series of the solution.
\end{proof}

\begin{theorem}
 If $J^x = v^x \id$, $J^y = v^y \id$, then the solution of \eqref{eq:generalsystem} with initial data $q(0,\vec x) = q_0(\vec x)$ is
\begin{align}
 q(t, \vec x) = q_0\Big(\vec x - \vec v_0\big(x - \vec v_0(\vec x)t\big)t\Big)  + \mathcal O(t^3)
\end{align}
where $\vec v_0(\vec x) = (v^x(0,\vec x), v^y(0,\vec x))$.
\end{theorem}
This is the same simple iteration as is in use for scalar conservation laws.
\begin{proof}
 We need to verify the presence of the term $\Delta t^2 \nabla q_i \Big|_{t=0} \cdot ((\vec v_0 \cdot \nabla) \vec v_0)$.
 \begin{align}
  \tilde q_i(t,\vec x) &:= q_{i,0}\Big(\vec x - \vec v_0\big(\vec x - \vec v_0(\vec x)t\big)t\Big) \\
  \del_t \tilde q_i(t,\vec x) &= \del_k  q_{i,0}\Big(\vec x - \vec v_0\big(\vec x - \vec v_0(\vec x)t\big)t\Big) \del_t \big(x_k - v^k_0(\vec x - \vec v_0(\vec x)t)t\big)\\
  &= - \del_k  q_{i,0}\Big(\vec x - \vec v_0\big(\vec x - \vec v_0(\vec x)t\big)t\Big)  \Big(- \del_\ell v^k_0\big(\vec x - \vec v_0(\vec x)t\big) v^\ell_0(\vec x) t + v^k_0\big(\vec x - \vec v_0(\vec x)t\big)\Big) \\
  \del_t^2 \tilde q_i(t,\vec x) \Big |_{t=0} &= 
  \del_\ell \del_k  q_{i,0}(\vec x)  v^\ell_0(\vec x) v^k_0(\vec x) +2 \del_k  q_{i,0}(\vec x)v^\ell_0(\vec x)   \del_\ell v^k_0(\vec x)  
 \end{align}
 Here $\del_k$ denotes $\del_0 \equiv \del_x$, $\del_1 \equiv \del_y$, and analogously for the components $v^k$ of $\vec v$. The first term is that of local linearization, and the second is precisely the extra term from Theorem \ref{thm:advectionapproxevoop}.
\end{proof}

In the Active Flux method therefore, instead of solving \eqref{eq:advectionfromeuler1}--\eqref{eq:advectionfromeuler2} via local linearization, two iterations are performed: One obtains first
\begin{align}
 \tilde {\vec x} := \vec x - \vec v_0\big(\vec x - \vec v_0(\vec x) t\big) t \label{eq:characteristicho}
\end{align}
and then, the reconstruction of $\rho, \vec v, p$ is evaluated at $\tilde {\vec x}$. Recall that this is a very special case, and does not lead to an increase in order in situations of non-commuting Jacobian matrices, i.e. it is not applicable to the acoustic operator. The difference of the present approach to that of \cite{roe18} is the absence of spatial derivatives. One can draw the following analogy to methods for ODEs: While the approach of \cite{maeng17,roe18}is similar to the Lax-Wendroff procedure converting higher time derivatives to space derivatives, the approach of \cite{barsukow19activeflux} as well as the prescription \eqref{eq:characteristicho} are more in the spirit of a Runge-Kutta method that generates the right terms in the Taylor series through nested evaluations.

\subsection{Bound preservation and limiting}

Since Active Flux is a high-order method, spurious oscillations need to be controlled. In general, for systems in multi-d there is no theory analogous to TVD for scalar conservation laws, and therefore there is no clear way to determine which oscillation is spurious and which is physical. For sure, however, undershoots that lead to negative values of density and/or pressure would lead to a loss of hyperbolicity, and thus definitely need to be prevented.

Here, the proposed approach is to use a posteriori limiting for point values and a priori limiting for the averages. Both strategies come with their respective advantages and disadvantages, and the choice of the present work tries to minimize the latter and keep the computational effort low. The a priori limiting of the averages is well-established (\cite{duan24,guermond16,kuzmin12,perthame96}) for a variety of methods. In principle, one could use a posteriori limiting (e.g. MOOD \cite{clain12,abgrall2023aftriangular}). To ensure conservation, the modified fluxes would require recomputation of a number of cells in the vicinity of the ``problematic'' one. This might result in an overall quicker algorithm, but would add algorithmic complexity. The a priori limiting for the cell averages is therefore preferred here (Section \ref{ssec:averageslimiting}).

For the point values (Section \ref{ssec:pointvalueslimiting}, a posteriori limiting is used: Since conservation is not an issue for point values, no recomputation of neighbours is needed and the approach is simple and very efficient. A priori limiting of point values was used in \cite{barsukow24afeuler}, but it was not possible to actually prove positivity preservation; in \cite{barsukow20swaf} a positivity preserving a priori limiting was developed in 1-d only.

Finally, to reduce other spurious oscillations, a shock sensor based procedure similar to that of \cite{duan24} is adopted (Section \ref{ssec:shocksensor}).

\subsubsection{Point value update} \label{ssec:pointvalueslimiting}

\begin{figure}
 \centering
 \includegraphics[width=0.3\textwidth]{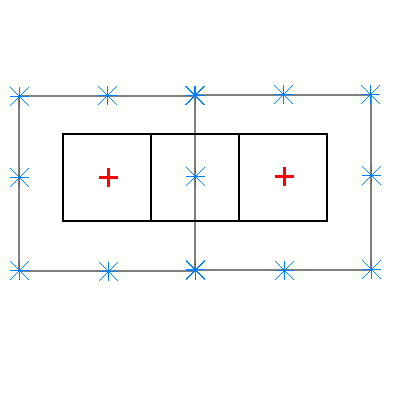} \hspace{2cm} \includegraphics[width=0.3\textwidth]{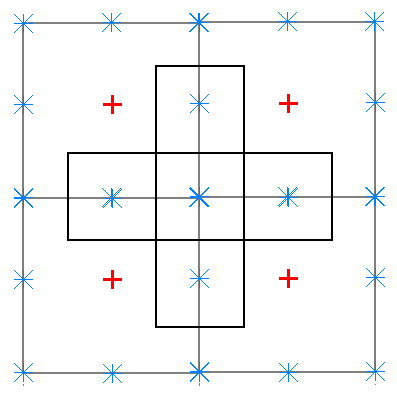}
 \caption{Control volumes on which the local Lax-Friedrichs scheme is applied for bound preservation. \emph{Left}: Update of the point value at an edge midpoint. \emph{Right}: Update of the point value at a node.}
 \label{fig:doflim}
\end{figure}

Since the new point values are required to obtain the fluxes for the average update, bound preservation for the point values is ensured first. Action is deemed necessary, if the new value is NaN (which can happen if the speed of sound has been evaluated at a location with negative density or pressure) or if the new value has a pressure or density that is lower then $\epsilon$. Here, $\epsilon = 10^{-10}$ is chosen. If limiting is needed, a local Lax-Friedrichs method based on the three values
\begin{align}
 &\bar q_{ij} && q_{i+\frac12,j} && \bar q_{i+1,j}
\end{align}
is used (and analogously for the other direction). This amounts to defining control volumes of size $\frac{\Delta x}{2} \times \frac{\Delta y}{2}$ centered at each degree of freedom (see Figure \ref{fig:doflim}), and identifying the average with the point value at cell center. 
Thus, the bound preserving point value at edge midpoint is
\begin{align}
 q_{i+\frac12,j}^{n+1,\text{lim}} = q_{i+\frac12,j}^n - \frac{\Delta t}{\Delta x/2} \left(\frac{f^x(\bar q^n_{i+1,j}) - f^x(\bar q^n_{ij})}{2} - \frac{\lambda^\text{HLL}}2 (\bar q^n_{i+1,j} - 2 q^n_{i+\frac12,j}  + \bar q^n_{ij})  \right)
\end{align}
 An analogous approach is used for the perpendicular edges, while for the node, four neighbours contribute: 
\begin{align}
 q_{i+\frac12,j+\frac12}^{n+1,\text{lim}} = q_{i+\frac12,j+\frac12}^n &- \frac{\Delta t}{\Delta x/2} \left(\frac{f^x(\bar q^n_{i+1,j+\frac12}) - f^x(\bar q^n_{i,j+\frac12})}{2} - \frac{\lambda^\text{HLL}}2 (\bar q^n_{i+1,j+\frac12} - 2 q^n_{i+\frac12,j+\frac12}  + \bar q^n_{i,j+\frac12})  \right) \\
 &\nonumber- \frac{\Delta t}{\Delta y/2} \left(\frac{f^y(\bar q^n_{i+\frac12,j+1}) - f^y(\bar q^n_{i+\frac12,j})}{2} - \frac{\lambda^\text{HLL}}2 (\bar q^n_{i+\frac12,j+1} - 2 q^n_{i+\frac12,j+\frac12}  + \bar q^n_{i+\frac12,j})  \right)
\end{align}

Since the point values are in primitive variables, prior to this calculation they need to be converted to conserved ones, and converted back afterwards. Observe that this is an \emph{a posteriori} limiting, and thus very inexpensive since it occurs rarely.

Finally, for all edges the states
\begin{align}
 q_{i+\frac12,j}^\text{HLL} = \frac{\bar q_{i+1,j} + \bar q_{ij}}{2} - \frac{1}{2\lambda^\text{HLL}_{i+\frac12,j} } \left(f^x(\bar q_{i+1,j}) - f^x(\bar q_{ij})\right)
\end{align}
and the fluxes
\begin{align}
 \hat f_{i+\frac12,j}^\text{HLL} = \frac{f^x(\bar q_{i+1,j}) + f^x(\bar q_{ij})}{2} - \frac{\lambda^\text{HLL}_{i+\frac12,j}}{2 } \left(\bar q_{i+1,j} - \bar q_{ij}\right) \label{eq:hllflux}
\end{align}
are computed (analogously for perpendicular edges). These states and fluxes correspond to considering only the grid of the averages, which has spacing $\Delta x$, and they are needed later for the bound preserving limiting of the averages.

As a practical simplification, the speed $\lambda^\text{HLL}_{i+\frac12,j}$ is computed using the point value at time $t^n$; a better way is described in \cite{duan24} and can be applied here as well. It is also taken as the maximum of the wave speed in both directions.

\subsubsection{Average update} \label{ssec:averageslimiting}

The procedure adopted here is the one from \cite{duan24}, which goes back to \cite{guermond16,kuzmin20,hajduk21} as well as \cite{perthame96,kuzmin12}: Given a high-order flux $\hat f_{i+\frac12,j}$ (from \eqref{eq:fluxquadraturex}) and a low-order bound-preserving flux $\hat f^\text{HLL}_{i+\frac12,j}$ (from \eqref{eq:hllflux}), a convex combination
\begin{align}
 \theta \hat f_{i+\frac12,j} + (1 - \theta) \hat f^\text{HLL}_{i+\frac12,j} &= \hat f^\text{HLL}_{i+\frac12,j} + \theta \Delta \hat f_{i+\frac12,j} \label{eq:convex}\\
 \Delta \hat f_{i+\frac12,j} &:= \hat f_{i+\frac12,j} - \hat f^\text{HLL}_{i+\frac12,j}
\end{align}
is found that is optimal, in the sense that it is as close to high-order as possible while maintaining bound-preservation. Since the method is 1-stage, there is no need to rely on an SSP property of the time integrator.

The condition for bound preservation is that the states
\begin{align}
 &q_{i+\frac12,j}^\text{HLL} \pm \theta \frac{\Delta \hat f_{i+\frac12,j}}{\lambda^\text{HLL}_{i+\frac12,j}} 
\end{align}
have positive density and pressure (see \cite{duan24} for a derivation). For the density, one thus imposes
\begin{align}
 \theta \Delta \hat f_{i+\frac12,j} = \begin{cases} \min(\Delta \hat f_{i+\frac12,j}, (\rho_{i+\frac12,j}^\text{HLL} - \epsilon)\lambda^\text{HLL}_{i+\frac12,j} ) & \Delta \hat f_{i+\frac12,j} > 0 \\
                                       \max (\Delta \hat f_{i+\frac12,j}, (\epsilon - \rho_{i+\frac12,j}^\text{HLL})\lambda^\text{HLL}_{i+\frac12,j} ) & \text{otherwise}
                                      \end{cases}
\end{align}
and computes the new flux of just the density according to \eqref{eq:convex}.

Next, the positivity of pressure leads to the inequalities
\begin{align}
 \theta^2 A \pm \theta B < C
\end{align}
which are ensured by enforcing (through $\theta \in [0,1]$)
\begin{align}
 \theta \big(\max(0, A) + |B|\big) < C
\end{align}
instead. The values of $A,B,C$ are
\begin{align}
 A &= \frac12 \left((\Delta f^{m_x})^2 + (\Delta f^{m_y})^2\right) - \Delta f^e \Delta f^\rho \\
 B &= \left(-{\vec m}^\text{HLL} \cdot \Delta f^{\vec m} + \Delta f^e \rho^\text{HLL} + \Delta f^\rho e^\text{HLL} - \Delta f^\rho \frac{\epsilon}{\gamma-1}\right)\lambda^\text{HLL}\\
 C &= \left ( - \rho^\text{HLL} \frac{\epsilon}{\gamma-1} - \frac12 \big( (m_x^\text{HLL})^2 + (m_y^\text{HLL})^2 \big) + e^\text{HLL} \rho^\text{HLL} \right) (\lambda^\text{HLL})^2
\end{align}
where, to save a tree, the indices have been dropped, and $\vec m$ denotes the momentum. With the $\theta$ thus obtained, \eqref{eq:convex} is used once again to compute the new flux of all the quantities.

Limiting with $\theta = 0$ is applied straight away if the high-order flux turns out to be NaN.

\subsubsection{Shock sensor based limiting} \label{ssec:shocksensor}

As shown in Section \ref{sec:numerical}, bound preservation alone leads to numerical solutions that display spurious oscillations. To control these, a procedure similar to that of \cite{duan24} is included as follows. First, two shock sensors are computed in every cell:
\begin{align}
 \phi_{ij,\text{Jameson}} &= \max\left( \frac{|p_{i+\frac12,j} - 2 p_{ij} + p_{i-\frac12,j}|}{p_{i+\frac12,j} + 2 p_{ij} + p_{i-\frac12,j}}, \frac{|p_{i,j+\frac12} - 2 p_{ij} + p_{i,j-\frac12}|}{p_{i,j+\frac12} + 2 p_{ij} + p_{i,j-\frac12}}\right )\\
 \phi_{ij,\text{div}} &= \max\left(0,  - \frac{\frac{u_{i+\frac12,j} - u_{i-\frac12,j}}{\Delta x} + \frac{v_{i,j+\frac12} -  v_{i,j-\frac12}}{\Delta y}}{\max\left(\frac{\max(|u_{i+\frac12,j}|, |u_{i-\frac12,j}|)}{\Delta x}, \frac{\max(|v_{i,j+\frac12}|, |v_{i,j-\frac12}|)}{\Delta y} \right )}\right)
\end{align}
$p_{ij}$ is the pressure obtained directly from the cell average.
Several modifications have been applied in comparison to the algorithm presented in \cite{duan24}. First, information local to the cell is used which results in a smaller stencil and a more concentrated shock sensor. Second, the divergence-based shock sensor (similar to that of Ducros) is normalized with an estimate of the jump of velocity across a distance of $\Delta x$. This is inspired by the fact that in strong shocks, the velocity jump is bounded by $\frac{\gamma+1}{\gamma-1}$ (just as the density); this factor is considered as order unity here. The modified normalization ensures that the sensor is not activated by noise; moreover its value has more physical relevance.

Then, after the modifications due to positivity of density and pressure have been applied, as described above, \eqref{eq:convex} is applied again on the resulting fluxes with its $\theta$ computed as
\begin{align}
 	\theta_{i+\frac12,j} = \exp\left (- \kappa \max\left(\phi_{i+1,j,\text{Jameson}}, \phi_{ij,\text{Jameson}}\right) 
							\max\left(\phi_{i+1j,\text{div}}, \phi_{ij,\text{div}}\right ) \right ) \in (0,1]
\end{align}
This follows closely \cite{duan24}, with the parameter value modified due to a different definition of the shock sensor. In all tests below, $\kappa = 10$ is used.

\section{Numerical results} \label{sec:numerical}

All numerical tests are performed with bound preservation turned on, and a CFL number of 0.45 with respect to $\mathrm{min}(\Delta x, \Delta y)$. Recall that the smallest separation of two degrees of freedom is half the edge length, and thus $\text{CFL} = 0.5$ is the stability limit (\cite{chudzik21}). For comparison, many of the test cases are those considered in \cite{barsukow24afeuler} for the semi-discrete Active Flux method. $\gamma = 1.4$ is used throughout.

\subsection{Convergence tests}

First, convergence tests are performed in order to assess the effect of the combination of a third-order space discretization, and a mixed-order time discretization.

\subsubsection{Contact wave}

\begin{figure}
 \centering
 \includegraphics[width=0.49\textwidth]{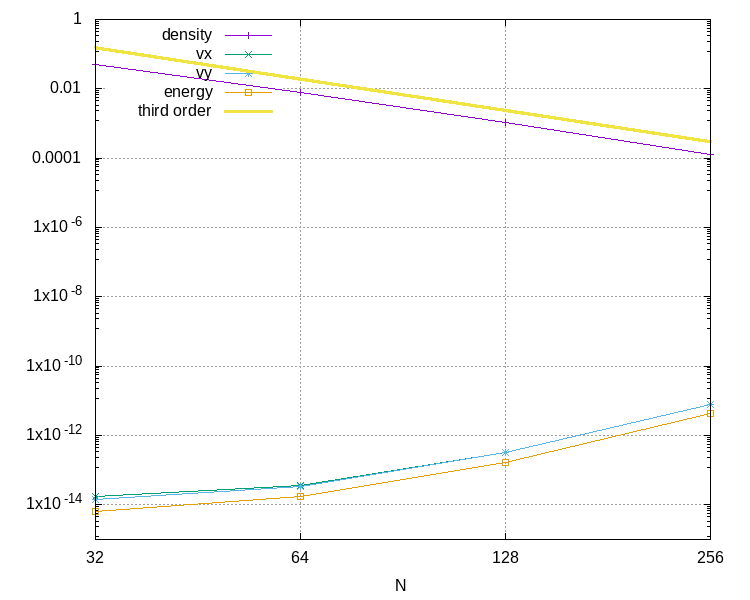} \hfill \includegraphics[width=0.49\textwidth]{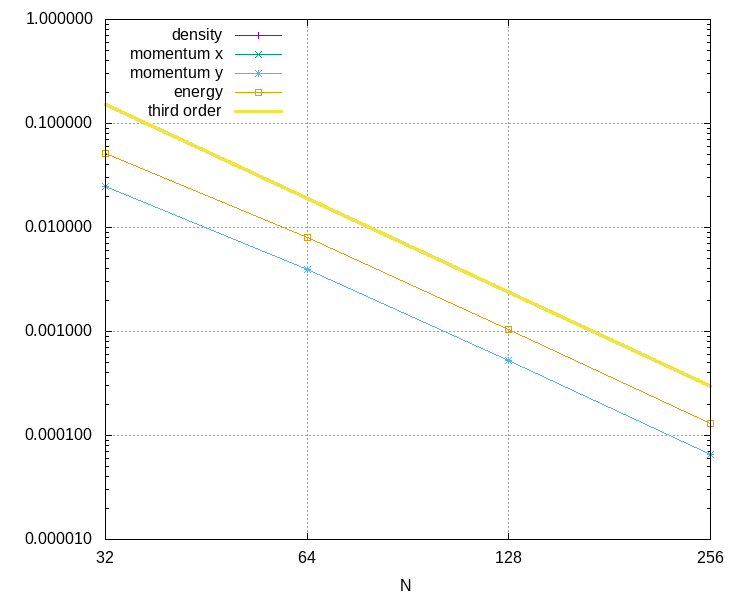} 
 \caption{Convergence study for the contact wave for different grid sizes $N \times N$. \emph{Left}: Errors of the point values on edges. \emph{Right}: Errors of the averages.}
 \label{fig:convergencehelzel}
\end{figure}

This test is from \cite{chudzik24}: The velocity and the pressure are constant:
\begin{align}
 \rho_0(x,y) &= \frac52 \exp\left( - 40 ((x-x_0)^2 + (y-y_0)^2) \right ) + \frac1{10} \\
 u_0(x, y) &= 1 \\
 v_0(x,y) &= 1\\
 p_0(x,y) &= 1
\end{align}
\cite{chudzik24} choose $x_0 = y_0 = -0.31$ on a domain $[-1,1]^2$. The exact solution is pure advection. Figure \ref{fig:convergencehelzel} shows the numerical error for various grid sizes at time $t=2$. One observes third order convergence (due to the high-order advection operator), and for the point values even a machine error for $p$ and $\vec v$. This is due to the usage of primitive variables, since both pressure and velocity are constant.

\subsubsection{Moving vortex}

\begin{figure}
 \centering
 \includegraphics[width=0.49\textwidth]{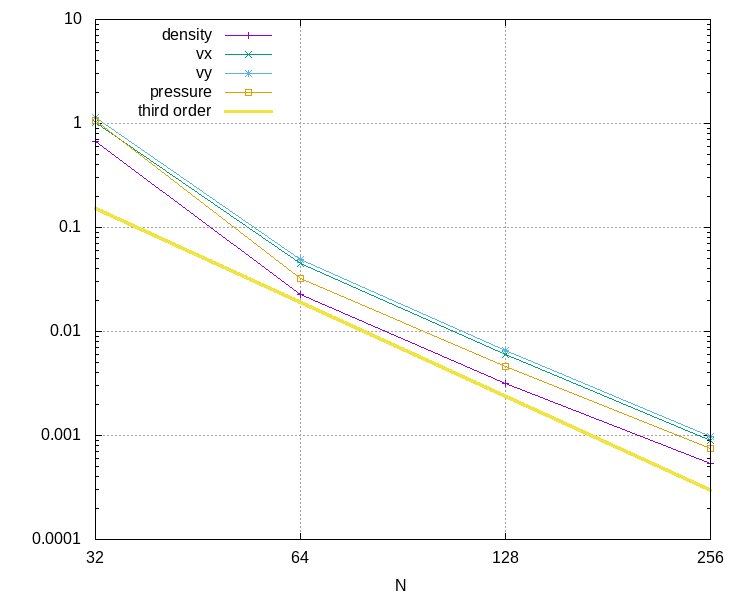} \hfill \includegraphics[width=0.49\textwidth]{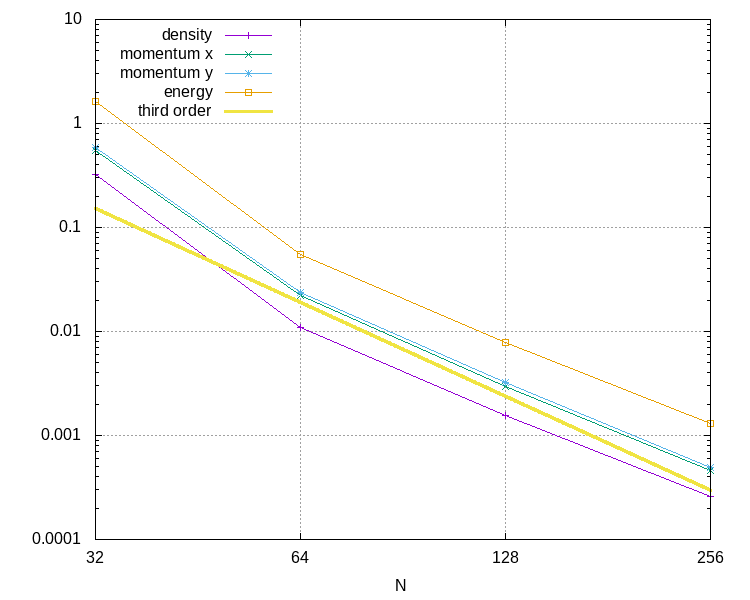} 
 \caption{Convergence study for the moving vortex for different grid sizes $N \times N$. \emph{Left}: Errors of the point values on edges. \emph{Right}: Errors of the averages.}
 \label{fig:convergence}
\end{figure}

A standard test case of an isentropic ($p = \rho^{\gamma}$) smooth moving vortex (\cite{hu99}, Example 3.3) is studied next. The initial data are
\begin{align}
\vecc{u_0(x, y)}{v_0(x, y)} &= \vecc{1}{1} + \frac{\Gamma }{2\pi}\exp\left(\frac{1 - r^2}2 \right) \vecc{-(y-5)}{x-5} \label{eq:convergence2u}\\
T_0(x,y) := \frac{p_0(x,y)}{\rho_0(x,y)} &= 1 - \frac{(\gamma-1) \Gamma^2}{8 \gamma \pi^2} \exp(1 - r^2) \label{eq:convergence2p}\\ 
r &:= \sqrt{(x - 5)^2 + (y -5)^2}
\end{align}
with the exact solution being advection with speed $(1, 1)$. However, this time, the pressure is not constant. $\Gamma = 5$ is used, as well as a periodic domain of $[0,10]^2$ and a final time of $t=10$, after which the vortex is back where it started from. The errors on different grids are shown in Figure \ref{fig:convergence}. One observes a convergence that is dominated by third order, which starts to deviate slightly towards finer grids. 

\subsection{Spherical shock tube}

\begin{figure}
 \centering
 \includegraphics[width=0.49\textwidth]{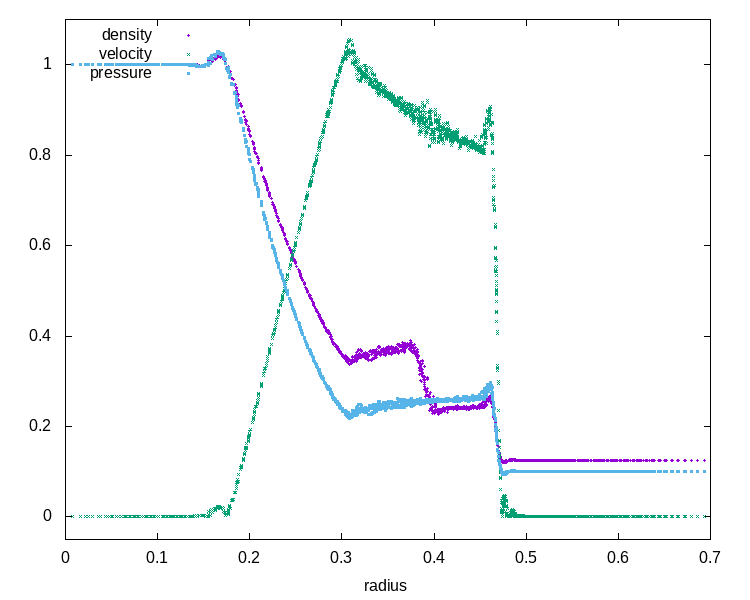} \hfill \includegraphics[width=0.49\textwidth]{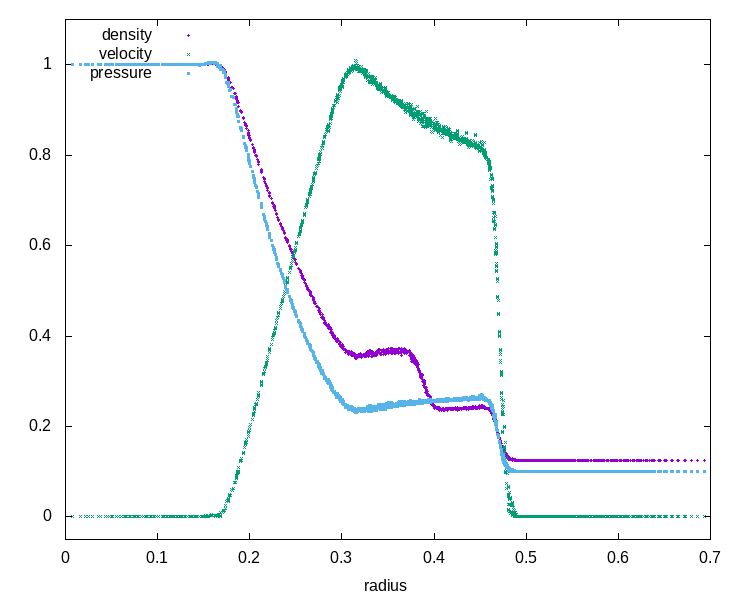}
 \caption{Scatter plots for the spherical Sod shock tube on a grid of $100\times 100$ cells. \emph{Left:} Bound preserving limiting used, but no limiting based on shock sensors. \emph{Right:} Same, but with shock sensor based limiting.}
 \label{fig:sod}
\end{figure}

Next, the performance of the method is assessed on flows with discontinuities. The first one is the spherical version of the Sod shock tube, i.e.
\begin{align}
 \rho_0(x, y) &= \begin{cases} 1 & r < 0.3 \\ 0.125 & \text{else} \end{cases} &
 p_0(x, y) &= \begin{cases} 1 & r < 0.3 \\ 0.1 & \text{else} \end{cases}  &
 \vec v_0(x, y) = 0
\end{align}
Without bound preserving limiting, this test fails; with it (but without the shock sensor based limiting), the procedure actually is activated only in the first two time steps (along the initial discontinuity). The results at time $t=0.1$ are shown in Figure \ref{fig:sod} (left). One observes some moderate oscillations and some scatter at the contact discontinuity. Part of it might be due to the imprint generated by the limiting at initial time. Upon using also the shock sensor based limiting (Figure \ref{fig:sod}, right), the oscillations are dramatically reduced.

\subsection{Multi-dimensional Riemann problems}

Next, some multi-dimensional Riemann problems from \cite{lax98} are used to assess the ability of the method to cope with more complex interactions between different wave families. The Riemann problems are set up such that outside the interaction region they result in a shock or a simple wave only. In fact, this is well-known to be a difficult situation for numerical methods. As is explained in Section 15.8.4 of \cite{leveque02}, whatever the numerical flux used to evolve a discontinuity between two values joined by, say, a Rankine-Hugoniot curve, the very first time step endows the cells in the vicinity with ($\Delta t$-dependent) values that are no longer on this curve in general. Thus, the \emph{second} time step creates waves of all families that propagate away as spurious artifacts. These waves are clearly visible in many simulation results of problems from \cite{lax98} in the literature, and they of course extend into the interaction region. Another reason for their presence is that the initial values given in \cite{lax98} fulfill the jump conditions only to a precision of about $10^{-3}$.

\begin{figure}
 \centering
 \includegraphics[width=0.49\textwidth]{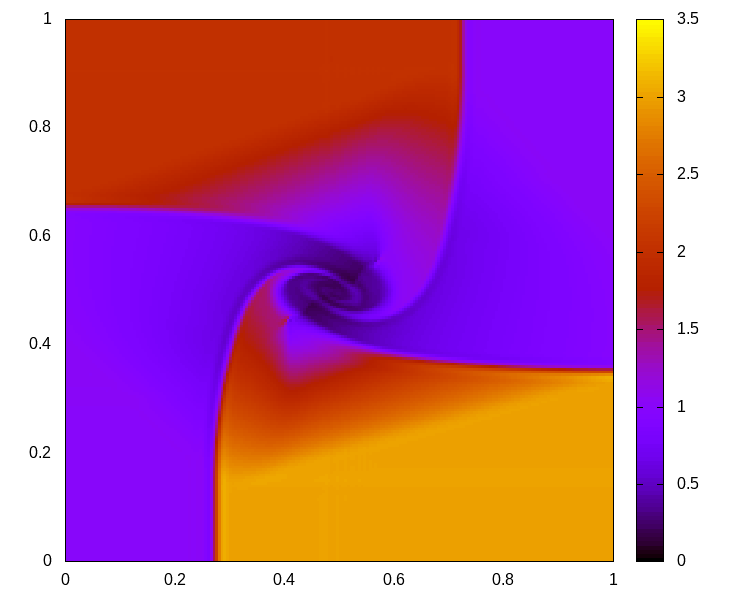} \hfill \includegraphics[width=0.49\textwidth]{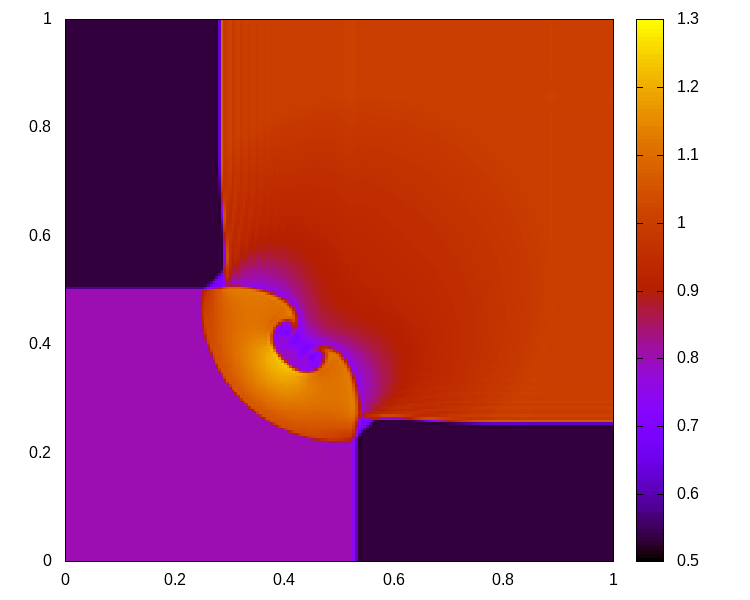} \\
 \includegraphics[width=0.49\textwidth]{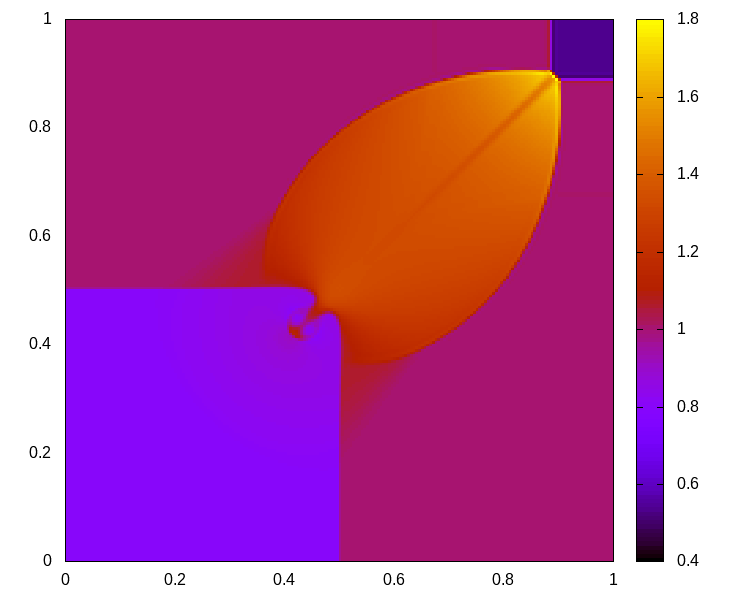} \hfill \includegraphics[width=0.49\textwidth]{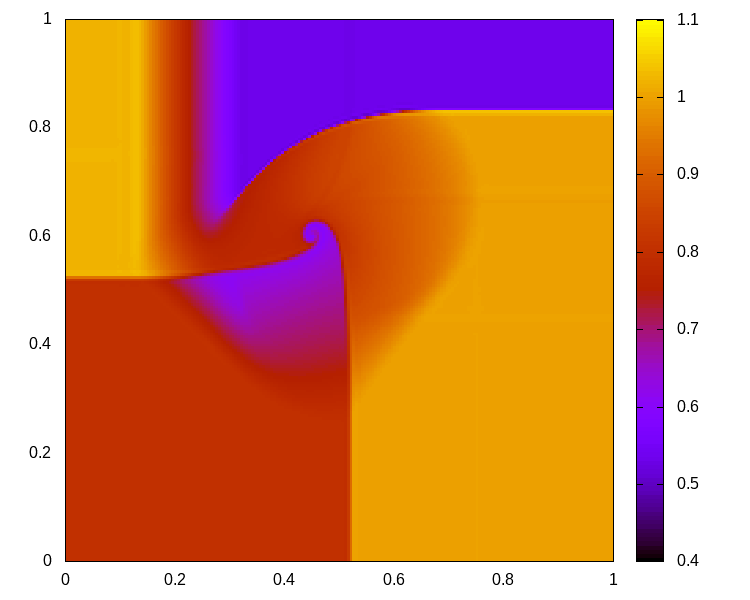}
 \caption{Several configurations from \cite{lax98} solved on a grid with $\Delta y = \Delta x = 1/200$. \emph{Left to right, top to bottom}: Configuration 6, $t=0.3$; Configuration  11,$t=0.3$; Configuration 12, $t=0.25$; Configuration 16 $t=0.2$. Color-coded is density.  Shock sensor based limiting is not used.}
 \label{fig:LL}
\end{figure}

\begin{figure}
 \centering
 \includegraphics[width=\textwidth]{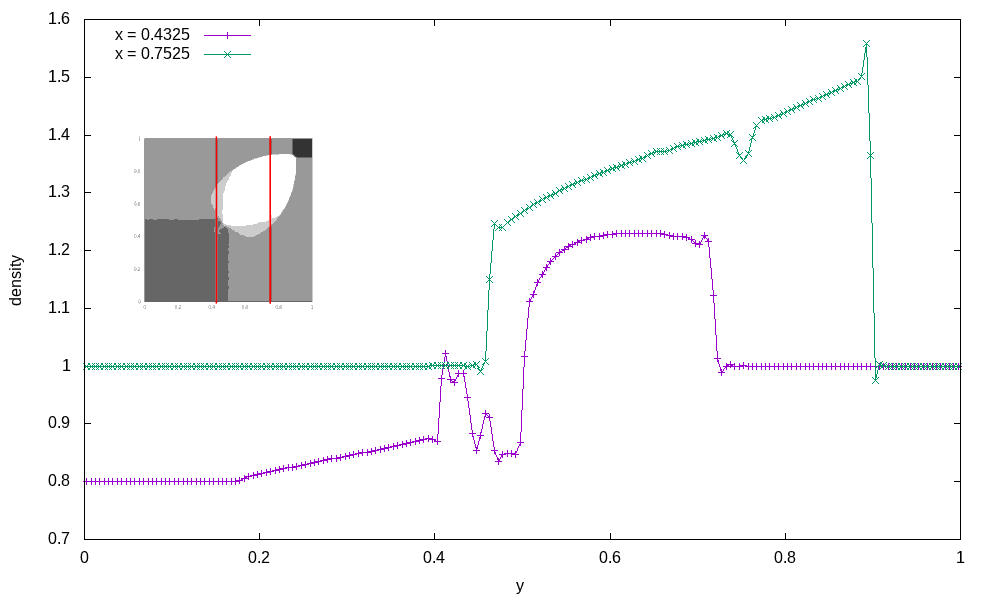}
 \caption{Cut through the solution of Configuration 12.}
 \label{fig:LL12cut}
\end{figure}

\begin{figure}
 \centering
 \includegraphics[width=0.33\textwidth]{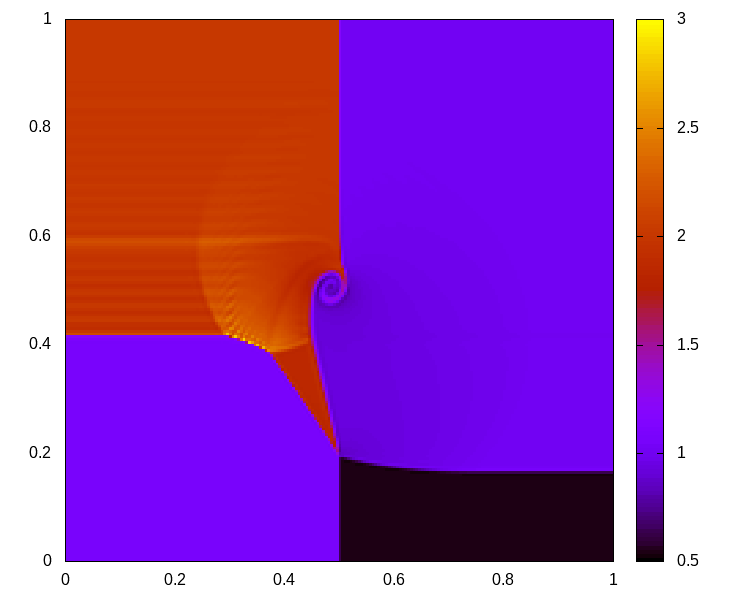}\includegraphics[width=0.33\textwidth]{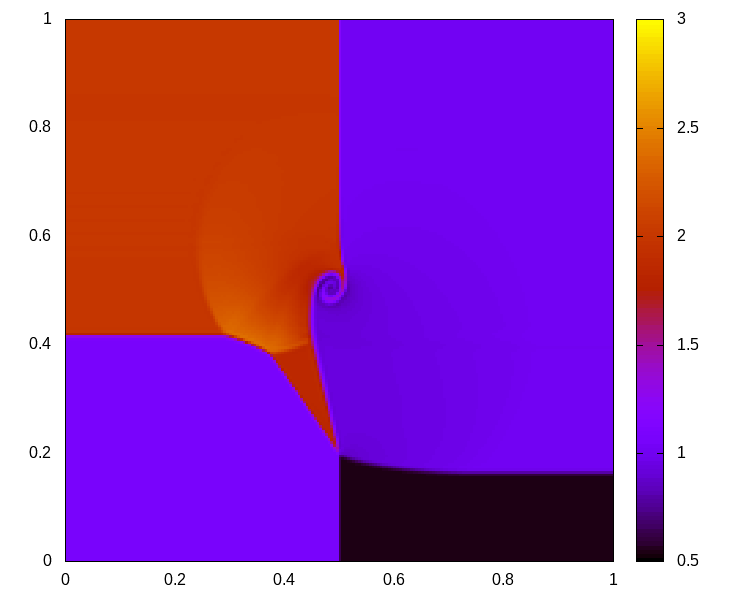}\includegraphics[width=0.33\textwidth]{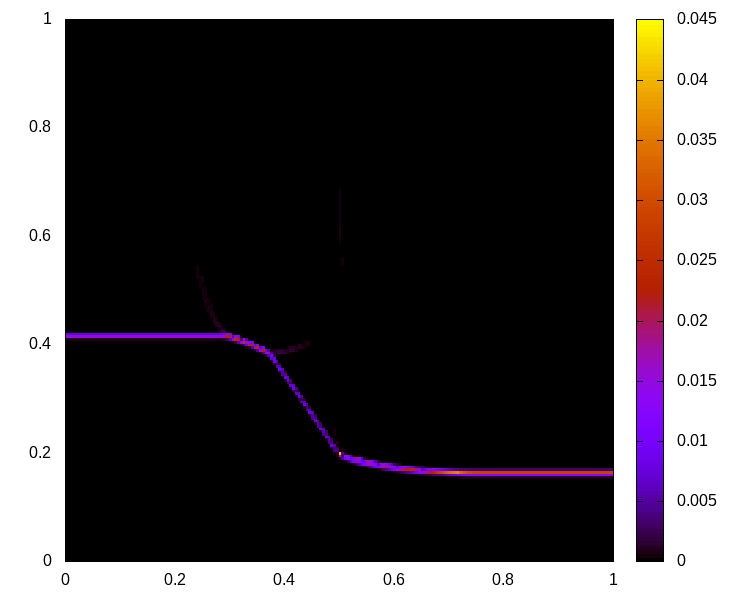}
\caption{Configuration 13 at time $t=0.3$. Color-coded is density. \emph{Left}: Bound preservation, but no shock sensor based limiting. Spurious oscillations propagate into the flow. \emph{Center}: Shock sensor based limiting used additionally; no oscillations present. \emph{Right}: Value $\phi_{\text{Jameson}} \phi_{\text{div}}$ of the product of the two shock sensors.}
 \label{fig:LL13}
\end{figure}

\begin{figure}
 \centering
 \includegraphics[width=0.49\textwidth]{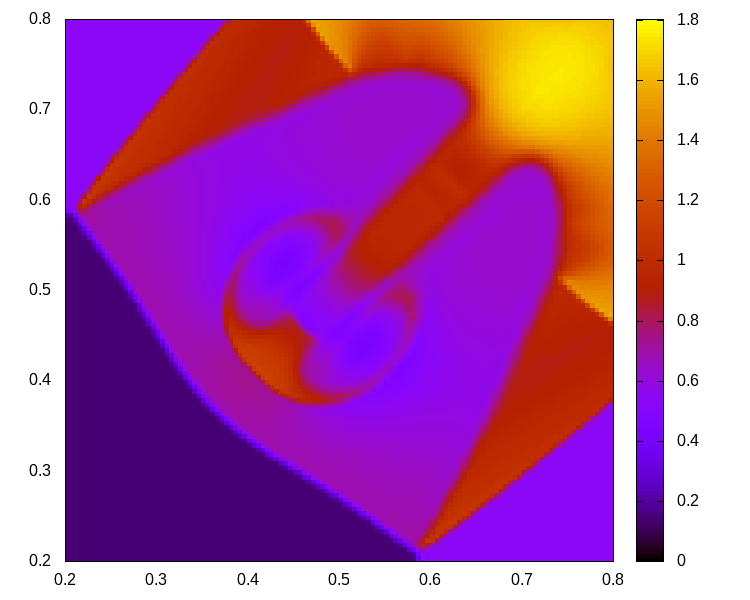} \hfill
 \includegraphics[width=0.49\textwidth]{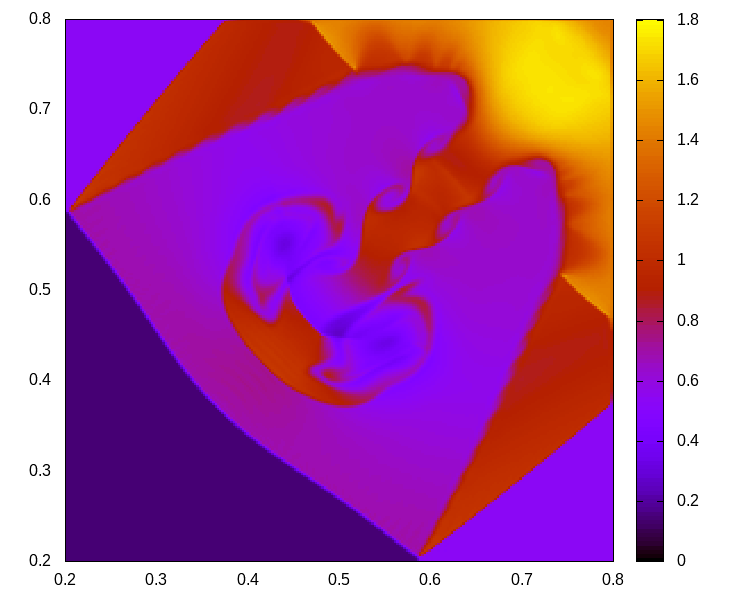}
 \caption{Modified configuration 3 from \cite{lax98} solved on a grid with $\Delta x = 1/200$ (\emph{left}) and $\Delta x = 1/400$ (\emph{right}) until $t=0.8$. Color-coded is density. Shock sensor based limiting is used.}
 \label{fig:LL3}
\end{figure}

A selection of multi-d Riemann problems is presented here (Figure \ref{fig:LL}). The number of grid cells is only half (in each direction) of what \cite{lax98} use. One observes that the features of the interaction region are resolved sharply. Figure \ref{fig:LL12cut} shows cuts through the simulation of Configuration 12 for two particular locations, in order to facilitate comparison with \cite{barsukow24afeuler}. In general, for the Riemann problems, no relevant differences to the results in \cite{barsukow24afeuler} can be established visually. Figure \ref{fig:LL13} demonstrates on Configuration 13 that spurious oscillations can propagate in the solution, and that the shock sensor based limiting is able to prevent them. Finally, Configuration 3, known for the development of Kelvin-Helmholtz-type instabilities at slip lines after long times is shown in Figure \ref{fig:LL3} for two different grid resolutions. This Configuration has been modified by adding 0.53 to both components of the initial velocity. The Galilei boost thus makes sure that the interaction region of the multi-dimensional Riemann Problem remains in the center of the computational domain.

\subsection{Low Mach number vortex}

\begin{figure}
 \centering
 \includegraphics[width=0.32\textwidth]{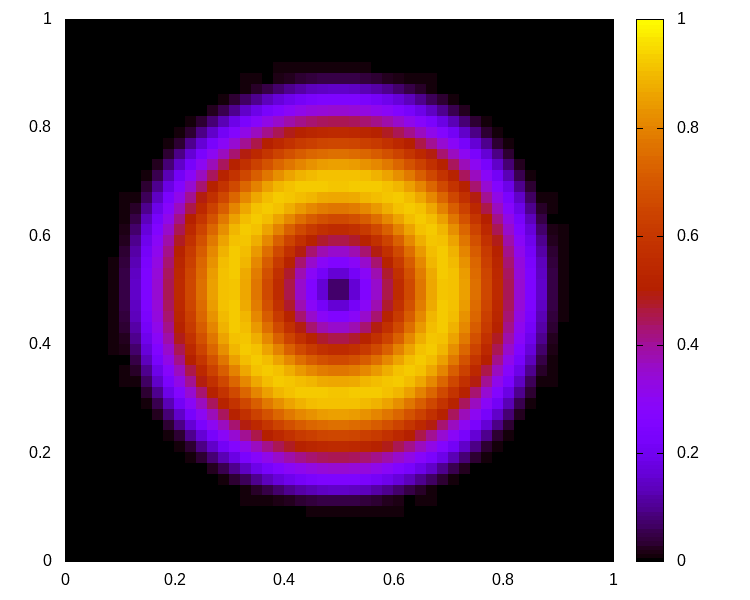} \hfill \includegraphics[width=0.32\textwidth]{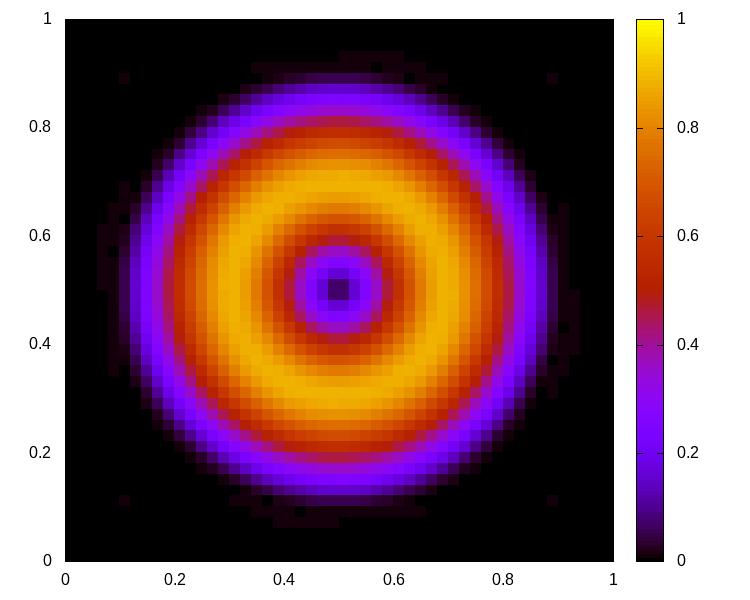}\hfill \includegraphics[width=0.32\textwidth]{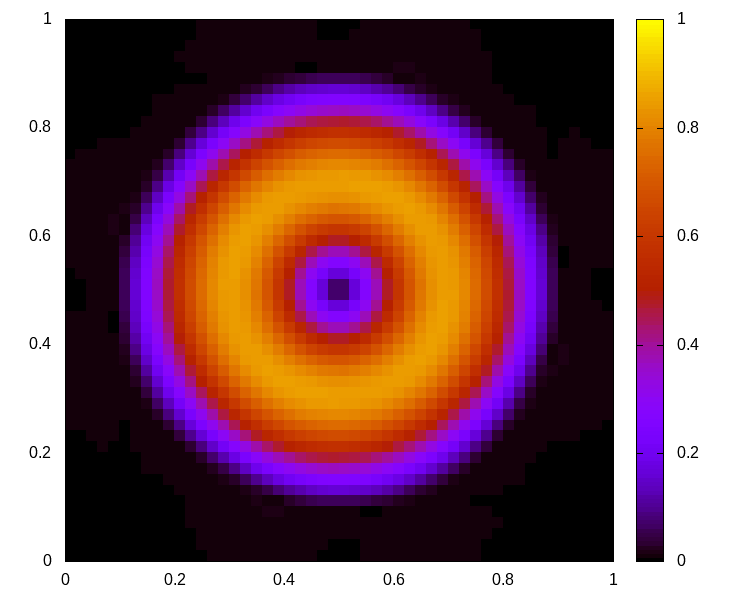}
 \\
 \includegraphics[width=0.32\textwidth]{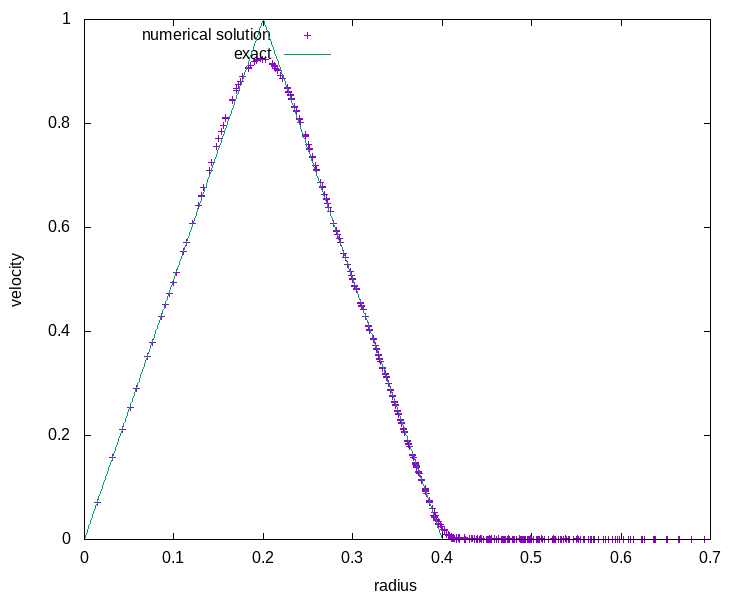} \hfill \includegraphics[width=0.32\textwidth]{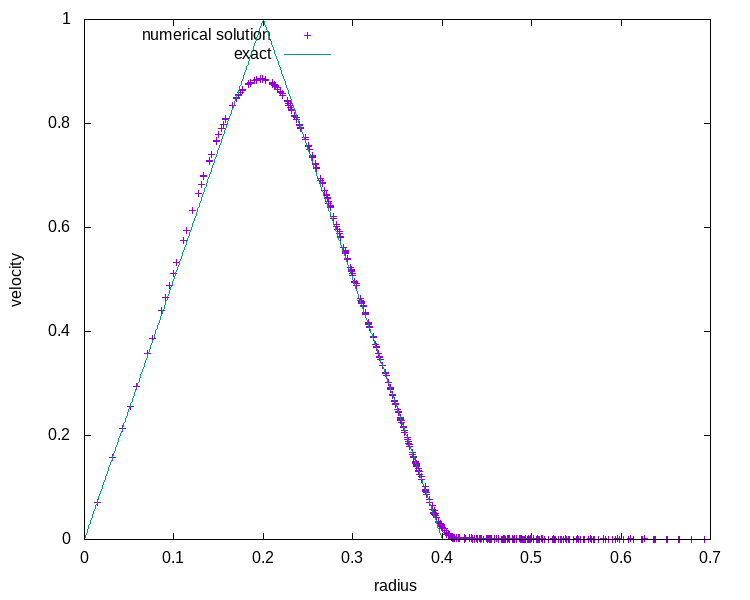}\hfill \includegraphics[width=0.32\textwidth]{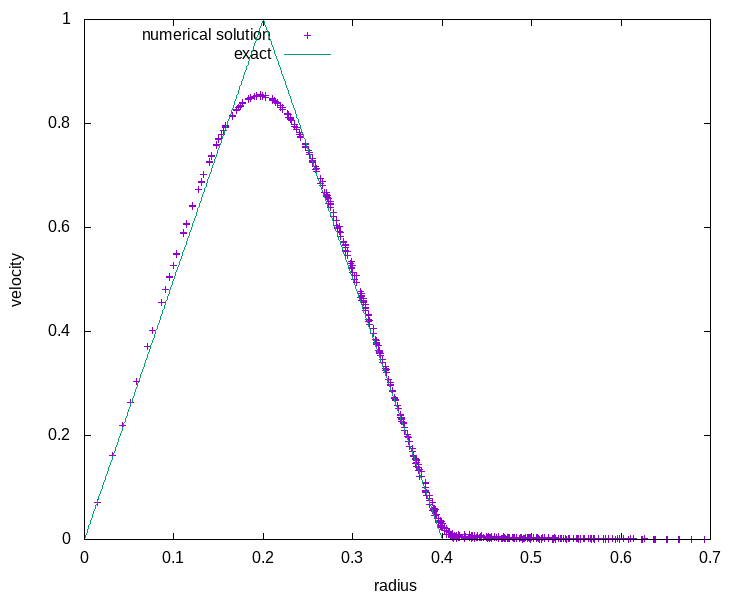}
 \caption{Numerical results for the stationary divergencefree Gresho vortex at different Mach numbers. \emph{From left to right}: $\epsilon = 10^{-1}, 10^{-2}, 10^{-3}$. \emph{Top}: Color coded is the magnitude of the velocity. \emph{Bottom}: Radial scatter plot of the magnitude of the velocity.}
 \label{fig:gresho}
\end{figure}

Low Mach number flows are of particular interest since numerical methods for compressible flow often require unreasonable grid refinement to resolve subsonic flow. A simple test case is that of a stationary, divergence-free vortex whose Mach number $\epsilon$ can be changed by modifying the background pressure, see \cite{barsukow16} for the details of the setup.
Results of a first-order Finite Volume method are available e.g. in \cite{barsukow20cgk}.
The simulation is evolved until $t=1$, which is when the quickest part of the flow has completed one revolution, while the sound wave crossing time is about $1/\epsilon$. Results for various choices of $\epsilon$ on a grid of $50 \times 50$ cells are shown in Figure \ref{fig:gresho}. The reader is referred to e.g. \cite{guillard99,dellacherie10,barsukow20cgk} for the background on the low Mach number problem. Non low-Mach-compliant methods dissipate the vortex on the sound wave crossing time, while low-Mach-compliant methods do so on an $\epsilon$-independent time scale in the limit $\epsilon \to 0$. This is observed in the experimental results: The diffusion (originating from the advection) is asymptotically independent of $\epsilon$. It is to be noted that this result is by no means trivial and that there is a large literature on ``fixing'' Finite Volume methods to achieve this. That Active Flux might be low-Mach-compliant has been observed in \cite{barsukow18activeflux} as a consequence of the stationarity-preserving property for linear acoustics, and low Mach compliance has been also experimentally observed in \cite{barsukow24afeuler} for the semi-discrete Active Flux method.

\subsection{Kelvin-Helmholtz instability}

The Kelvin-Helmholtz instability allows to qualitatively assess the numerical diffusion of a method in subsonic conditions, but in a significantly more complex flow than for an isolated, stationary vortex. Methods that are not low-Mach-compliant are usually artificially stabilizing the setup, preventing the formation of any vortices. The relevant aspect of the next tests is to assess the mere fact of vortex formation. If the initial vorticity is concentrated on a line, then the evolution of the vortices cannot be reasonably studied in the inviscid setting anyway. However, the results can still be compared across different methods in order to qualitatively assess the amount of numerical diffusion. This is the purpose of the test in Section \ref{ssec:sound} where the instability is triggered by the passage of a sound wave. The test in Section \ref{ssec:khhd} has smooth initial vorticity, and therefore until a certain time a well-defined time evolution to which numerical methods converge upon mesh refinement (\cite{leidi24}).

\subsubsection{Smooth shear flow} \label{ssec:khhd}

The initial data are as follows (\cite{leidi24}):
\begin{align}
\rho_0(x, y) &:= \gamma + R  (1 - 2\eta(y)) & p_0(x,y) &:= 1\\
v_{y,0}(x,y) &:= \delta M \sin(2 \pi x) &
v_{x,0} &:= M (1 - 2 \eta(y))
\end{align}
with
    		\begin{align}
\eta(y) &:=
		\begin{cases}
			\frac12 \left(1 + \sin\left(16\pi \left(y + \frac14\right)\right)\right)
		 & -\frac{9}{32} \leq y < -\frac{7}{32} \\
			1
		&  -\frac{7}{32}\leq y < \frac{7}{32} \\
			\frac12 \left(1 - \sin\left(16 \pi \left(y - \frac14\right)\right)\right)
		&  \frac{7}{32} \leq y < \frac{9}{32} \\
		0	
		& \text{else}
    \end{cases}
\end{align}

\begin{figure}
 \centering
 \includegraphics[width=0.7\textwidth]{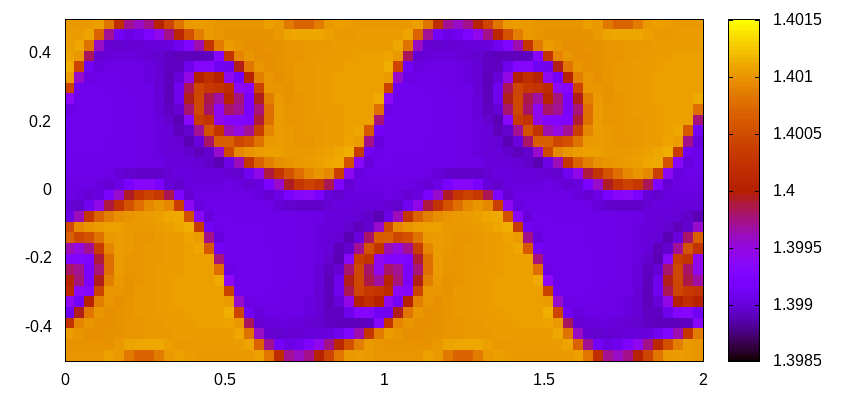} \\
 \includegraphics[width=0.7\textwidth]{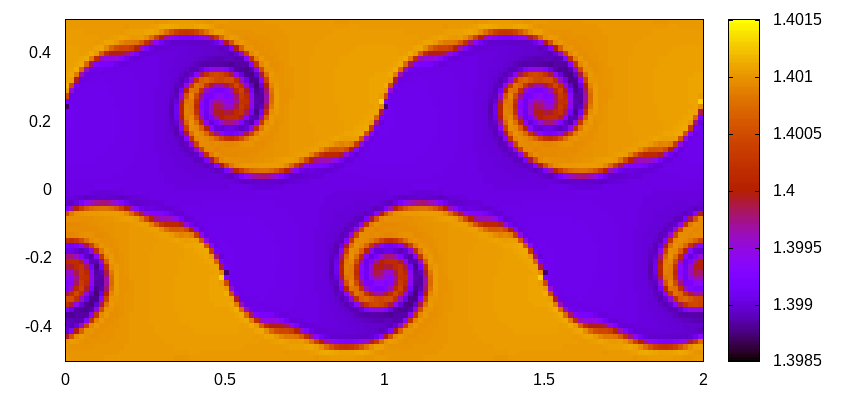} \\
 \includegraphics[width=0.7\textwidth]{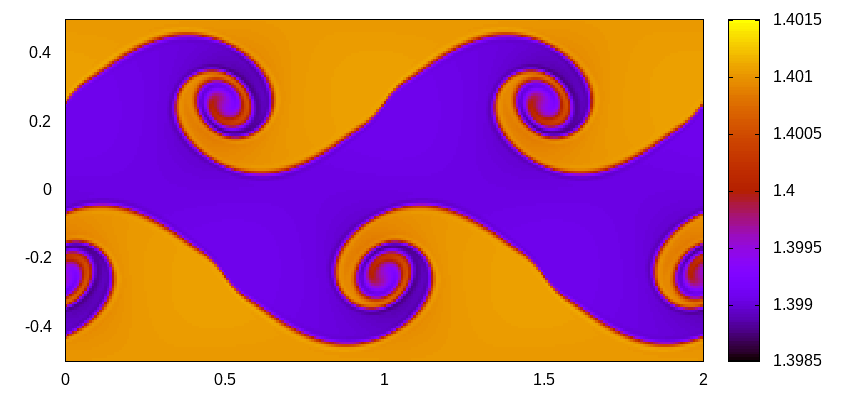}
 \caption{Smooth shear flows on different grids. \emph{Top to bottom}: $64 \times 32$, $128 \times 64$, $256 \times 128$}
 \label{fig:khhd}
\end{figure}

Here, $R = 10^{-3}$ and $\delta = 0.1$. $M$ is a parameter governing the Mach number of the flow, chosen here as $M = 0.01$. The setup is considered on a grid covering $[0,2] \times [-\frac12,\frac12]$ subject to periodic boundaries until a final time of $0.8/M$. Figure \ref{fig:khhd} shows the results on several different grids. One observes that the new method is able to resolve well the development of the vortices (see \cite{leidi24} for a comparison with other methods). For low resolution, some spurious features on the interface away from the vortices are visible, which are common among many methods and sometimes even manifest themselves as secondary, spurious vortices (\cite{brown95}). This is not the case here and the features disappear upon mesh refinement.

\subsubsection{Instability triggered by the passage of a sound wave} \label{ssec:sound}

Another type of a Kelvin-Helmholtz instability is a test case proposed in \cite{munz03}, as a way to assess the ability of the numerical method to cope with acoustic waves in a nearly incompressible flow. A saw-tooth profile of density is hit by an acoustic wave, which induces a shear flow. Since the initial discontinuity in the density translates into a discontinuity in tangential speed, this setup displays ever more features upon grid refinement. Thus, the actual evolution of the vortices cannot be considered as relevant, however, the time scales at which they appear and their minimal size give qualitative information on the numerical diffusion of the method in the subsonic regime. The initial data are given and explained in detail in \cite{barsukow18activeflux}. The reader is also referred to this work for comparison of the results with those of the semi-discrete Active Flux, as well as to the original work \cite{munz03} (for comparison with Finite Volume methods) and \cite{peraire11} (for comparison with a DG method).

\begin{figure}
 \centering
 \includegraphics[width=\textwidth]{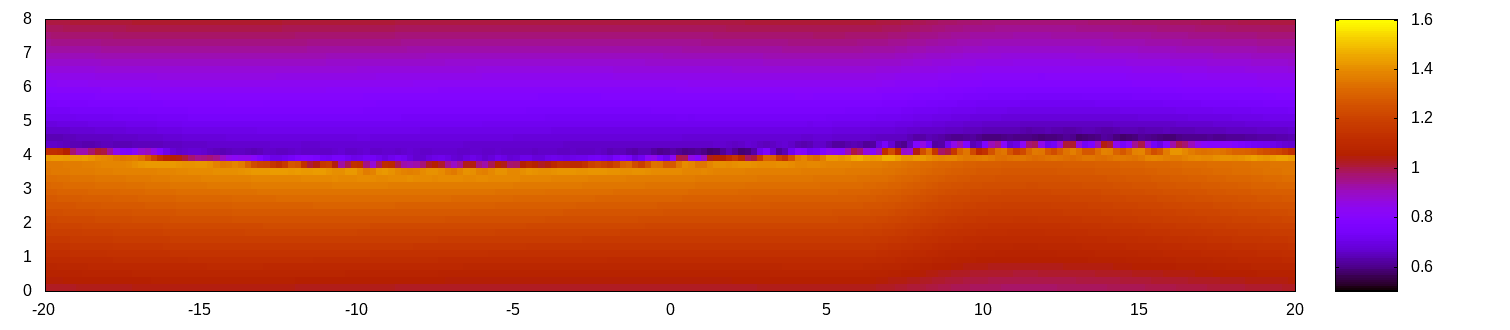} \\
 \includegraphics[width=\textwidth]{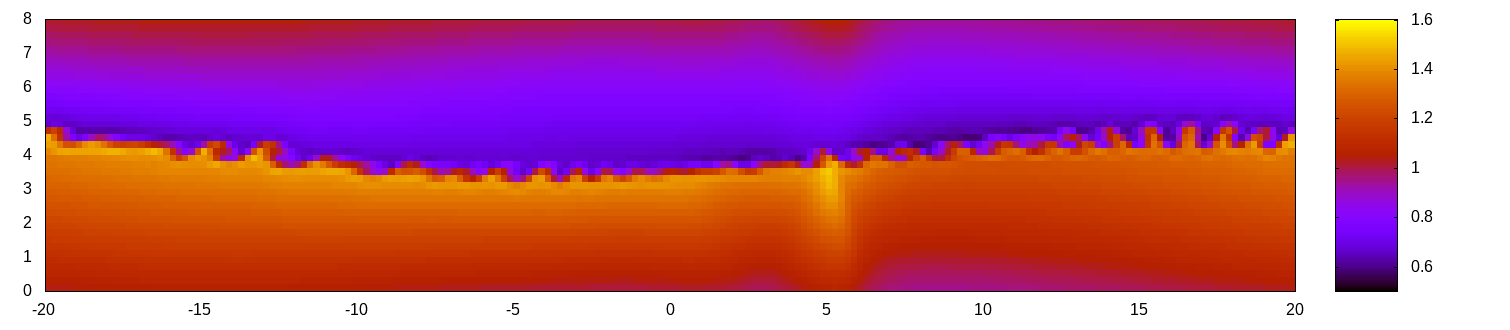} \\
 \includegraphics[width=\textwidth]{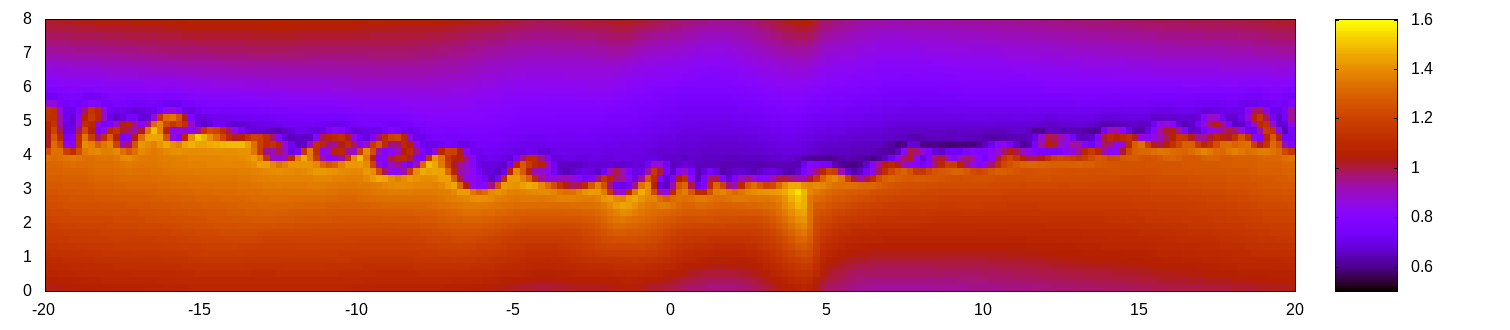} \\
 \includegraphics[width=\textwidth]{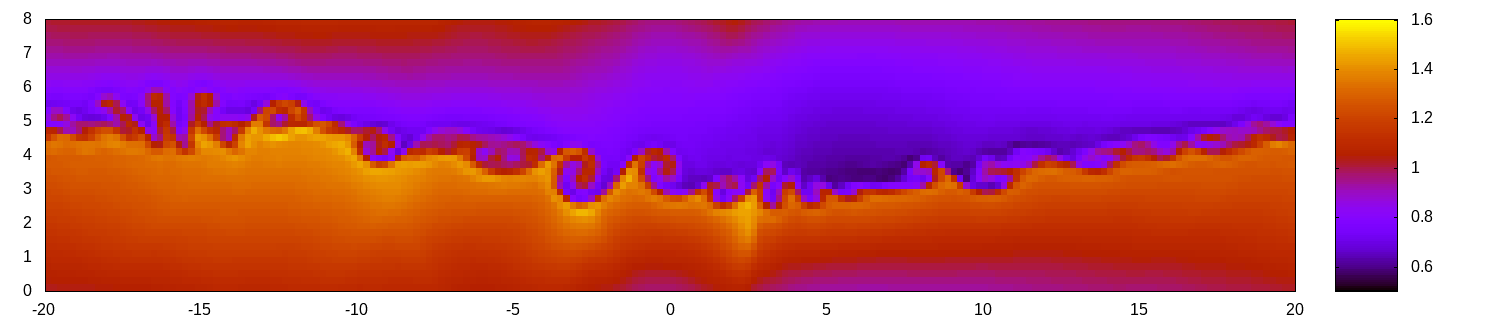}\\
 \includegraphics[width=\textwidth]{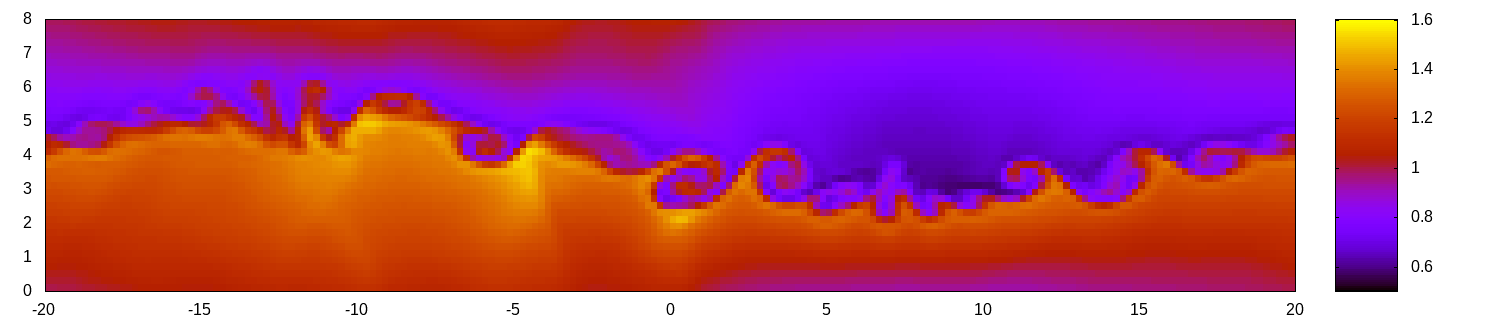}
 \caption{Kelvin-Helmholtz instability triggered by the passage of a sound wave, on a grid of $200 \times 40$ cells. Density is shown at times $t=3,6,9,12,15$.}
 \label{fig:khsoundwave}
\end{figure}

\begin{figure}
 \centering
 \includegraphics[width=\textwidth]{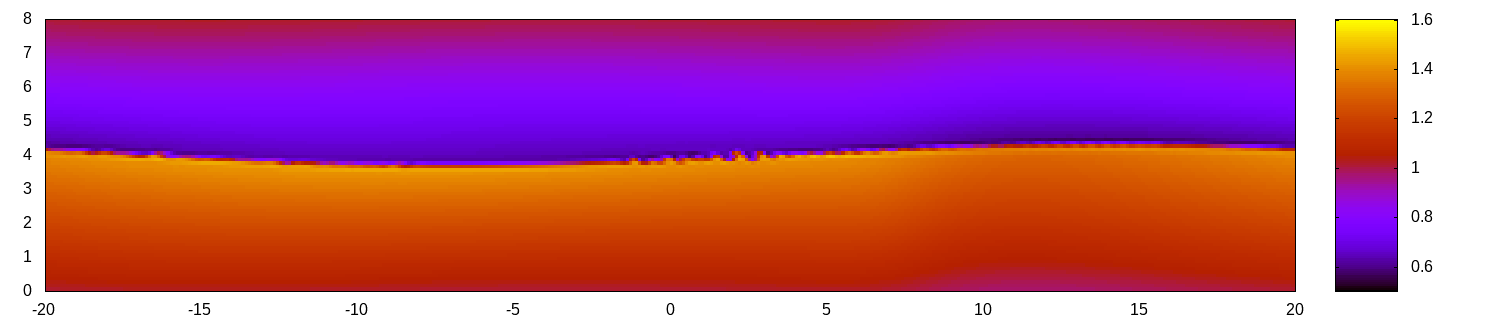} \\
 \includegraphics[width=\textwidth]{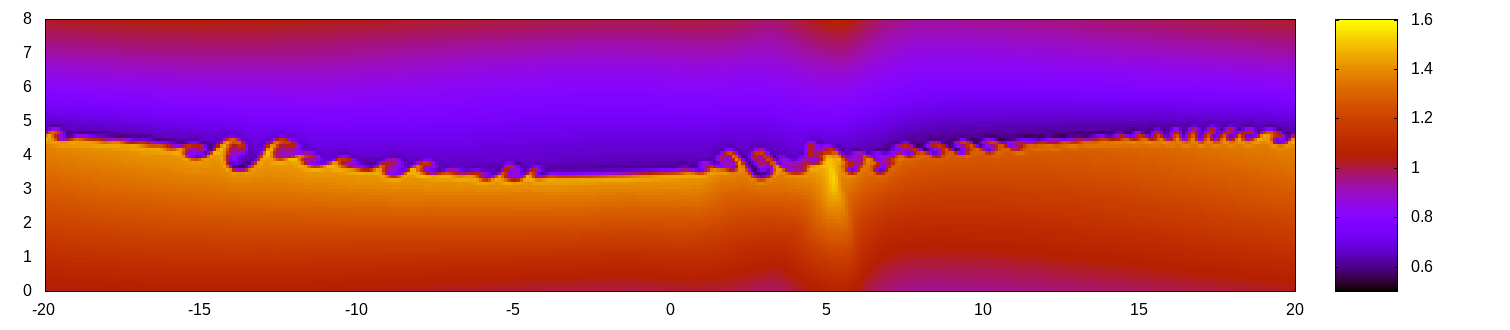} \\
 \includegraphics[width=\textwidth]{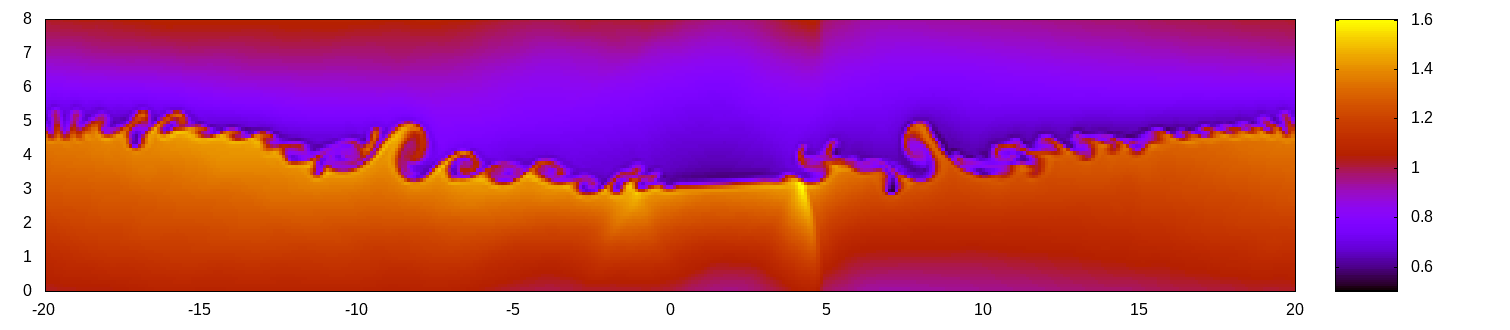} \\
 \includegraphics[width=\textwidth]{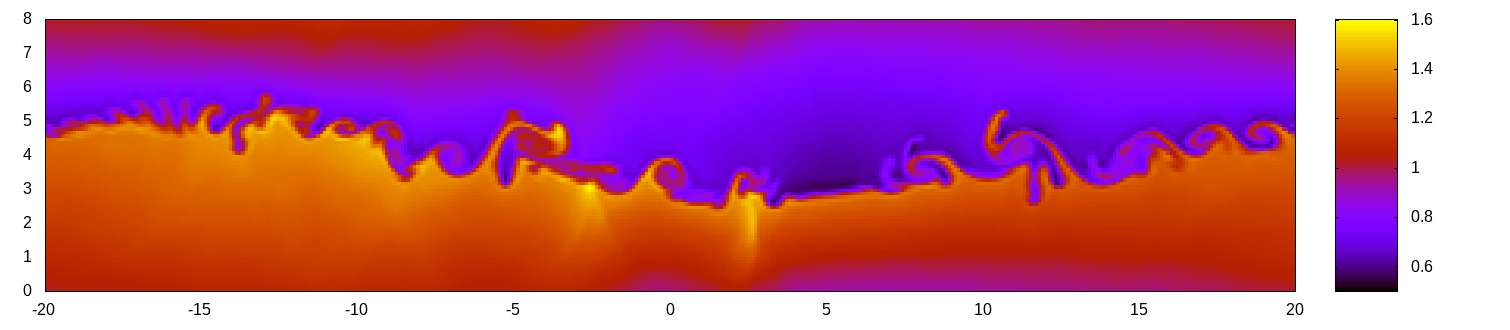}\\
 \includegraphics[width=\textwidth]{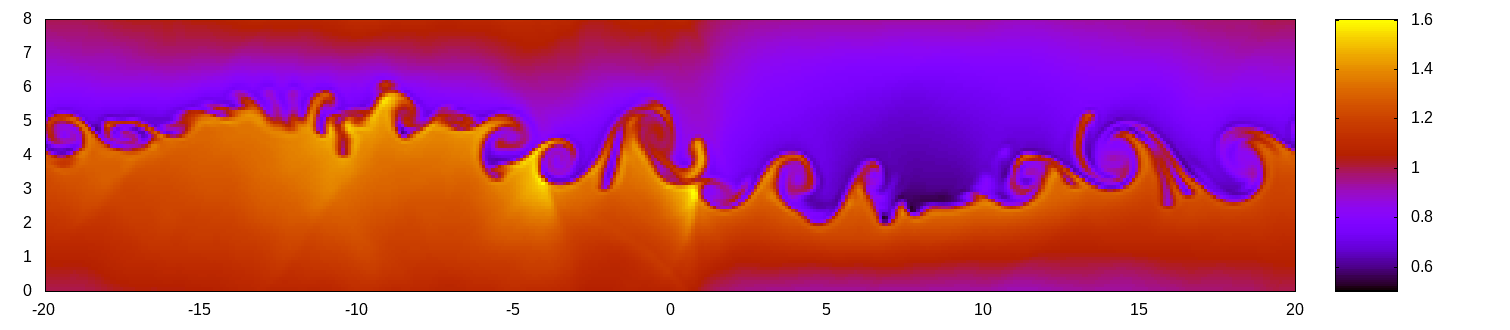}
 \caption{Same as Figure \ref{fig:khsoundwave}, but on a grid of $400 \times 80$ cells.}
 \label{fig:khsoundwave2}
\end{figure}

Figure \ref{fig:khsoundwave} shows the results for a rather coarse grid of $200 \times 40$ cells. One observes that the vortices of the instability appear much earlier for the present method than for the semi-discrete method from \cite{barsukow24afeuler}, and more homogeneously across the entire interface. This effect appears even stronger if one takes into account the finer resolved grids used in \cite{barsukow24afeuler}. It thus seems that the present method has significantly less artificial diffusion than its semi-discrete counterpart. Figure \ref{fig:khsoundwave2} shows a simulation on $400 \times 80$ cells.

\begin{figure}
 \centering
 \includegraphics[width=\textwidth]{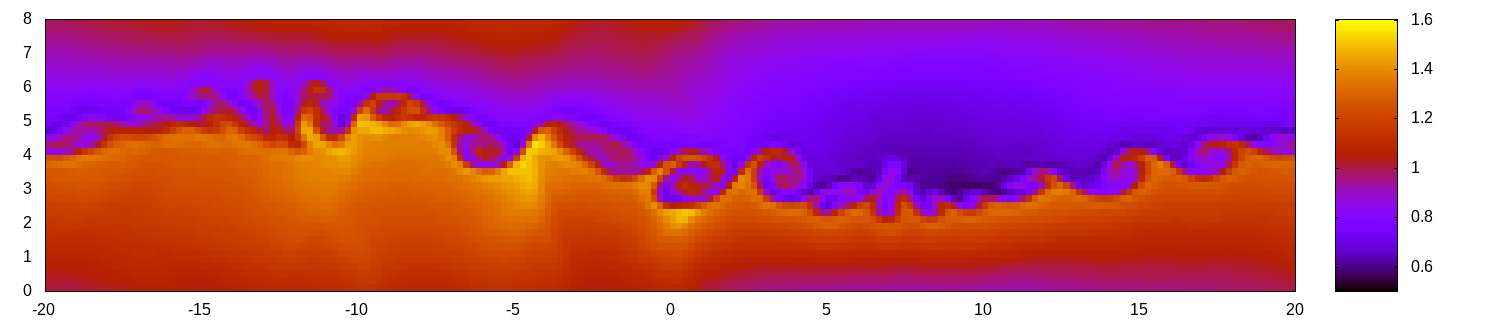}\
 \caption{Kelvin-Helmholtz instability triggered by the passage of a sound wave, on a grid of $200 \times 40$ cells, upon usage of the low-order advection operator \eqref{eq:characteristicsimple}. Density is shown at time $t=15$.}
 \label{fig:khsoundwavelo}
\end{figure}

Figure \ref{fig:khsoundwavelo} shows the results of a simulation that uses the simple, second-order evolution operator \eqref{eq:characteristicsimple} for advection, instead of the higher-order version \eqref{eq:characteristicho} used otherwise. There are visible differences (for the vortex at $x\simeq12$, for example), yet the overall behaviour is similar.

\section{Conclusions and outlook}

Attempts to develop multi-dimensional Active Flux methods based on evolution operators very quickly reached success for the linear acoustic equations. For the significantly more complicated case of the multi-dimensional Euler equations, several suggestions exist. The main decisions to take in their design, and hence the main differences between existing approaches are whether to use splitting into advection and acoustics (yes), whether to use evolution operators that are exact for the corresponding linearized problems (yes) and whether to use operators that all are sufficiently high-order accurate in time (no) -- the terms in brackets indicating the decisions taken in the present work. In practice, the latter question needs to be answered by taking into account the (probably) significant computational effort associated to a fully third-order evolution operator, and the ability to answer it in the affirmative depends on whether it is at all possible to find this operator. Unfortunately, it does not seem to be currently the case, even though a roadmap was traced in \cite{roe18}. While these considerations apply to the multi-d Euler equations, and also the acoustic suboperator, they do not apply to the advective sub-operator, which has a much simpler nature. As is shown in the present work it can easily be solved to the required order of accuracy. Concerning its accuracy in space, the method is of third order.

The usage of operators that are exact for linearized problems guarantees von Neumann stability of the method. Moreover, since the low Mach compliance for the Euler equations is linked to stationarity preservation for acoustics, a link between the good results of the present method at low Mach number and the corresponding analysis of stationarity preservation in \cite{barsukow18activeflux} can be made immediately. The extensive details of an efficient implementation of the exact evolution operator for linear acoustics presented here for the first time hopefully can serve as a reference in future.

Besides the aspects already mentioned, the current work features an example of point values and the reconstruction in primitive variables, as well as a blending of an a priori and a posteriori bound preservation. The usage of primitive variables simplifies the splitting and the corresponding evolution operators greatly, and in certain situations can result in enhanced accuracy.

Numerical results demonstrate that the new method is at least as able to resolve multi-dimensional Riemann problems as the recently developed semi-discrete Active Flux method (\cite{barsukow24afeuler}), while it might have less diffusion at least in the subsonic regime. It has a higher maximal CFL.

Future work will be devoted to further comparison between the semi-discrete Active Flux and the one based on evolution operators. Other ways to achieve higher order and the possible associated benefits will be investigated.

\bibliographystyle{alpha}
\newcommand{\etalchar}[1]{$^{#1}$}

\appendix

\section{Spherical means} \label{sec:sphericalmean}

The acoustic evolution operator involves spherical means of the initial data, on which certain differential operators act. 
Apart from the usual spherical mean of a scalar function, for vector-valued functions there exist additional spherical means of interest that involve the unit normal $\vec n = \veccc{\cos \phi \sin \theta}{\sin \phi \sin \theta}{\cos \theta}$. For integration over the angular coordinates in three spatial dimensions, $\dd  \Omega :=  \dd \phi \dd \theta \sin \theta$ is used occasionally as an abbreviation.

\subsection{General definitions}

\begin{definition}
 Consider a scalar-valued function $f \colon \mathbb R^3 \to \mathbb R$ and a vector valued function $\vec v \colon \mathbb R^3 \to \mathbb R^3$, $\vec v = (v_1, \ldots, v_n)$. 

 \begin{enumerate}[i.]
  \item 
The spherical mean $M[f] \colon \mathbb R^3 \times \mathbb R^+_0 \to \mathbb R$ of $f$ is 
\begin{align}
 M[f](\vec x,r) &:= \frac{1}{4\pi} \int_{S^2} \dd \vec y \, f(\vec x + r \vec y) 
\end{align}
\item 
The spherical mean $M[\vec v] \colon \mathbb R^3 \times \mathbb R^+_0 \to \mathbb R^3$ of a vector-valued function $\vec v$, is defined componentwise:
 \begin{align}
  M[\vec v](\vec x , r) := \left( M[v_1](\vec x , r), \ldots, M[v_n](\vec x , r) \right)
 \end{align}
 In components, this can be expressed as $M[\vec v]_i = M[v_i]$.
 \item The spherical mean $M[\vec n f] \colon \mathbb R^3 \times \mathbb R^+_0 \to \mathbb R^3$ of $f$ is 
 \begin{align}
  M[\vec n f](\vec x , r) &:= \frac{1}{4\pi} \int_{S^2} \dd \vec y \, \vec y f(\vec x + r \vec y)
 \end{align}
 \item The spherical mean $M[\vec n \cdot \vec v] \colon \mathbb R^3 \times \mathbb R^+_0 \to \mathbb R$ of $\vec v$ is 
 \begin{align}
  M[\vec n \cdot \vec v](\vec x , r) := \frac{1}{4\pi} \int_{S^2} \dd \vec y \, \vec y \cdot \vec v(\vec x + r \vec y)
 \end{align}
 \item The spherical mean $M[\vec n (\vec n \cdot \vec v)] \colon \mathbb R^3 \times \mathbb R^+_0 \to \mathbb R^3$ of $\vec v$ is 
 \begin{align}
  M[\vec n(\vec n \cdot \vec v)](\vec x , r) := \frac{1}{4\pi} \int_{S^2} \dd \vec y \, \vec y (\vec y \cdot \vec v(\vec x + r \vec y))
 \end{align}
\end{enumerate}
 \end{definition}

 Observe that the normalization is such that $M[1] = 1$. 
For actual computations, the spherical mean is best expressed in spherical coordinates, e.g.
\begin{align}
 M[f](\vec x, r) &= \frac{1}{4\pi} \int_{0}^{2\pi} \dd \phi \int_0^\pi \dd \theta \, \sin \theta \, f(x + r \cos \phi \sin \theta, y + r \sin \phi \sin \theta, z + r \cos \theta) \\
 M[\vec n f](\vec x, r) &= \frac{1}{4\pi} \int_{0}^{2\pi} \dd \phi \int_0^\pi \dd \theta \, \veccc{\cos \phi \sin \theta}{\sin \phi \sin \theta}{\cos \theta} \sin \theta \, f(x + r \cos \phi \sin \theta, y + r \sin \phi \sin \theta, z + r \cos \theta)
\end{align}
etc. The definition for a spherical mean over a smaller angular domain is analogous. 

\subsection{Spherical means in two spatial dimensions}

On two-dimensional domains, integration over subsets of the full domain of $\phi$ are of interest (because the data might be piecewise defined). 

\begin{definition}
The angular domain $W = [\phi_\text{min}, \phi_\text{max}) \times [0, \pi]$ is called a \emph{wedge} $W \subset S^2$. Spherical means whose $(\phi, \theta)$ integration domain is $W$ are denoted by $M^W$:
\begin{align}
 M^W[f](\vec x,r) &:= \frac{1}{2(\phi_\text{max}-\phi_\text{min})} \int_{\phi_\text{min}}^{\phi_\text{max}} \dd \phi \int_0^\pi \dd \theta \, \sin \theta \, f(x + r \cos \phi \sin \theta, y + r \sin \phi \sin \theta)
\end{align}
with $f : \mathbb R^2 \to \mathbb R$. An analogous definition holds for other types of spherical means such as $M^W[\vec n f]$.
\end{definition}

\begin{remark}
Observe that it only makes sense to define $M^W[f]$ for $f$ that do not depend on $z$.
\end{remark}

\begin{remark}
 In two spatial dimensions it is natural to only consider the first two components of $\vec n$. This will not be made explicit in the notation. This is consistent, since for any $f$ that does not depend on $z$, 
 \begin{align}
  M^W[\vec n f]_z &= \frac{1}{2(\phi_\text{max}-\phi_\text{min})} \int_{\phi_\text{min}}^{\phi_\text{max}} \dd \phi \int_0^\pi \dd \theta \, \sin \theta \cos \theta \, f(x + r \cos \phi \sin \theta, y + r \sin \phi \sin \theta) \\
  &\overset{c := \cos \theta}{=} \frac{1}{2(\phi_\text{max}-\phi_\text{min})} \int_{\phi_\text{min}}^{\phi_\text{max}} \dd \phi \int_{-1}^1 \dd c \, c \, f(x + r \cos \phi \sqrt{1 - c^2}, y + r \sin \phi \sqrt{1 - c^2}) \
  &= 0
 \end{align}
 because an antisymmetric function is integrated over a symmetric interval.
\end{remark}

The normalization $M^W[1] = 1$ is easily confirmed:
\begin{align}
 \int_{W \subset S^2} \dd \vec y = 
 \int_{\phi_\text{min}}^{\phi_\text{max}} \dd \phi \int_0^\pi \dd \theta \, \sin \theta = 
 (\phi_\text{max} - \phi_\text{min}) \int_{-1}^1 \dd(\cos \theta) = 2 (\phi_\text{max} - \phi_\text{min})
\end{align}

\subsection{Spherical means for polynomial data in two space dimensions}

In principle, it is possible to evaluate everything analytically. However, in particular for angular domains that are not the entire sphere, the expressions become complicated. Here, some formulae are given that allow to (pre)compute the expressions exactly for polynomial initial data in 2 spatial dimensions.

For polynomial initial data, the key object to compute therefore is the spherical mean
\begin{align}
 M^W[n_x^\alpha n_y^\beta x^a y^b]
\end{align}

\subsubsection{Integration over $\theta$}

In the following we compute the spherical mean of a polynomial function in two spatial dimensions. For polynomial functions, the dependence on $\theta$ will be of the form $\sin^a \theta$, which justifies the following definition.

\begin{definition}
Denote by $\displaystyle \eta_a := \int_0^\pi \dd \theta \sin^{a} \theta$.
\end{definition}

\begin{theorem}[\cite{gradshteyn63}, 3.621 (3. and 4.)]
One has 
\begin{align}
 \eta_a = \sqrt{\pi} \frac{\Gamma \left ( \frac{a+1}{2}\right )}{\Gamma \left( \frac{a+2}{2} \right )} = \frac{(a-1)!!}{a!!} \cdot \begin{cases} \pi & \text{$a > 0$ even}\\ 2 & \text{$a > 0$ odd} \end{cases}
\end{align}
where the double factorial is defined as
\begin{align}
 A!! := \begin{cases} A (A-2)(A-4) \cdots 4 \cdot 2 & \text{$A$ even} \\ A (A-2)(A-4) \cdots 3 \cdot 1 & \text{$A$ odd} \end{cases}
\end{align}
\end{theorem}

This means that $
 (\eta_a)_{a=0,1,\ldots} = \pi, 2, \frac{\pi}{2}, \frac43, \frac{3\pi}{8}, \frac{16}{15}, \frac{5\pi}{16}, \frac{32}{35}, \ldots
$. Below, we will need the special case 
\begin{align}
 \eta_{i+j+1}= \frac{(i+j)!!}{(i+j+1)!!} \cdot \begin{cases} \pi & \text{$i+j+1$ even}\\ 2 & \text{else} \end{cases}
\end{align}

\subsubsection{Integration over $\phi$}

\begin{definition}
 Denote by $\displaystyle \mu^W_{ij} := \int_{\phi_\text{min}}^{\phi_\text{max}}  \dd \phi \cos^{i} \phi \sin^{j} \phi $.
\end{definition}

\begin{theorem}
 We have 
 \begin{align}
  \mu^W_{ij}&= \frac{1}{2^{i+j} } \sum_{m=0}^i \sum_{n=0}^j {i \choose m} {j \choose n} (-1)^{n+j} \\&  \phantom{mmmm}\times \begin{cases} (-1)^{j/2} (\phi_\text{max} - \phi_\text{min}) & i+j =2 (m + n), j \text{ even}\\ &\\ 
 \displaystyle \left.  \frac{(-1)^{j/2} \sin( [ i+j -2 m - 2 n ] \phi)}{  i+j -2 m - 2 n } \right |_{\phi_\text{min}}^{\phi_\text{max}} & i+j \neq2 (m +  n), j \text{ even} \\ 
 \displaystyle \left.  \frac{(-1)^{\frac{j-1}{2} } \cos( [ i+j -2 m - 2 n ] \phi)}{  i+j -2 m - 2 n } \right |_{\phi_\text{min}}^{\phi_\text{max}} & i+j \neq2 (m +  n), j \text{ odd} \\ 0 & \text{else}\end{cases}
 \end{align}
with the notation $f(\phi) \Big|_{\phi_\text{min}}^{\phi_\text{max}} = f(\phi_\text{max}) - f(\phi_\text{min})$.
\end{theorem}
\begin{proof}
The integral can be evaluated to
\begin{align}
\mu_{ij}^W = \int_{\phi_\text{min}}^{\phi_\text{max}}  \dd \phi  \cos^{i} \phi   \sin^{j} \phi&= \int_{\phi_\text{min}}^{\phi_\text{max}}  \dd \phi  \left(  \frac{\ee^{\ii \phi} - \ee^{- \ii \phi}}{2 \ii}  \right )^{j}  \left(  \frac{\ee^{\ii \phi} + \ee^{- \ii \phi}}{2}  \right )^{i}\\
 &= \frac{1}{2^{i+j} \ii^j} \int_{\phi_\text{min}}^{\phi_\text{max}}  \dd \phi \,\, \ee^{\ii (i+j) \phi}  (  1 - \ee^{- 2\ii \phi}  )^{j}  (  1 + \ee^{-2 \ii \phi } )^{i}\\
 &=  \frac{1}{2^{i+j} \ii^j} \sum_{m=0}^i \sum_{n=0}^j {i \choose m} {j \choose n} (-1)^{n} \int_{\phi_\text{min}}^{\phi_\text{max}}  \dd \phi \,\,  \ee^{\ii [ i+j -2 m - 2 n ] \phi} 
\end{align}
Now,
\begin{align}
 \int_{\phi_\text{min}}^{\phi_\text{max}}  \dd \phi \,\,  \ee^{\ii [ i+j -2 m - 2 n ] \phi} = \begin{cases} \phi_\text{max} - \phi_\text{min} & i+j -2 m - 2 n = 0\\ &\\ \left. \displaystyle \frac{\ee^{\ii [ i+j -2 m - 2 n ] \phi}}{\ii [ i+j -2 m - 2 n ]} \right |_{\phi_\text{min}}^{\phi_\text{max}} & \text{else} \end{cases} \end{align}
and thus
\begin{align}
 \mu_{ij}^W &= \frac{1}{2^{i+j} (-1)^j} \sum_{m=0}^i \sum_{n=0}^j {i \choose m} {j \choose n} (-1)^{n} \begin{cases} \ii^j (\phi_\text{max} - \phi_\text{min}) & i+j -2 m - 2 n = 0\\ &\\ \left. \displaystyle -\frac{\ii^{j+1} \ee^{\ii [ i+j -2 m - 2 n ] \phi}}{  i+j -2 m - 2 n } \right |_{\phi_\text{min}}^{\phi_\text{max}} & \text{else} \end{cases}
\end{align}
As the integral is real, we only need to include real terms into the sum:
\begin{align}
  \mu_{ij}^W &= \frac{1}{2^{i+j} } \sum_{m=0}^i \sum_{n=0}^j {i \choose m} {j \choose n} (-1)^{n+j} \\&  \phantom{mmmm}\times \begin{cases} (-1)^{j/2} (\phi_\text{max} - \phi_\text{min}) & i+j =2 (m + n), j \text{ even}\\ &\\ 
 \displaystyle \left.  \frac{(-1)^{j/2} \sin( [ i+j -2 m - 2 n ] \phi)}{  i+j -2 m - 2 n } \right |_{\phi_\text{min}}^{\phi_\text{max}} & i+j \neq2 (m +  n), j \text{ even} \\ 
 \displaystyle \left.  \frac{(-1)^{\frac{j-1}{2} } \cos( [ i+j -2 m - 2 n ] \phi)}{  i+j -2 m - 2 n } \right |_{\phi_\text{min}}^{\phi_\text{max}} & i+j \neq2 (m +  n), j \text{ odd} \\ 0 & \text{else}\end{cases}
 \end{align}
 
\end{proof}

The first few values of $\mu_{ij}^W$ are:

\begin{center}
 \begin{tabular}{c|c|l} 
 $i$ & $j$ & $\mu_{ij}^W \cdot \Delta \phi$ \\\hline\hline
0 & 0 & $\Delta\phi$\\
0 & 1 & $-\Delta\cos(\phi)$\\
0 & 2 & $\frac1{4} \left(- \Delta\sin(2\phi)+2 \Delta\phi\right)$\\ \hline
1 & 0 & $\Delta\sin(\phi)$\\
1 & 1 & $-\frac1{4} \Delta\cos(2\phi)$\\
1 & 2 & $\frac1{4} \left(- \frac{\Delta\sin(3\phi)}{3}+\Delta\sin(\phi)\right)$\\\hline
2 & 0 & $\frac1{4} \left(\Delta\sin(2\phi)+2 \Delta\phi\right)$\\
2 & 1 & $-\frac1{4} \left(\frac{\Delta\cos(3\phi)}{3}+ \Delta\cos(\phi)\right)$\\
2 & 2 & $\frac1{16} \left(- \frac{\Delta\sin(4\phi)}{2} +2 \Delta\phi  \right)$
\end{tabular}
\end{center}

Here, $\Delta \phi = \phi_\text{max} - \phi_\text{min}$ and $\Delta\sin(n\phi) = \sin(n \phi_\text{max}) - \sin(n \phi_\text{min})$ etc.

\subsubsection{Assembling the spherical mean}

By linearity and by the assumption that we are concerned with the two-dimensional case only it is sufficient to compute the spherical mean at $\vec x = (x,y,0)$ of
\begin{align}
 n_x^\alpha n_y^\beta x^a y^b
\end{align}
Therefore it is necessary to evaluate the following spherical mean
\begin{align}
 M^W[n_x^\alpha n_y^\beta x^ay^b](\vec x, r) &=
 \frac{1}{2 (\phi_\text{max}-\phi_\text{min})} \\\nonumber &\!\!\!\!\!\!\!\!\!\!\!\!\!\!\!\!\!\!\!\!\!\!\!\!\!\!\!\!\!\!\!\!\!\!\!\!\!  \cdot \int_{\phi_\text{min}}^{\phi_\text{max}} \dd \phi \int_0^\pi \dd \theta \,\,\sin \theta (\cos \phi \sin \theta)^\alpha (\sin\phi\sin\theta)^\beta (x + r \cos \phi \sin \theta)^a (y + r \sin \phi \sin \theta)^b \\
 &\!\!\!\!\!\!\!\!\!\!\!\!\!\!\!\!\!\!\!\!\!\!\!\!\!\!\!\!\!\!\!\!\!\!\!\!\! = \frac{1}{2 (\phi_\text{max}-\phi_\text{min})} \int_{\phi_\text{min}}^{\phi_\text{max}} \dd \phi \int_0^\pi \dd \theta \,\, \sin^{1+\alpha+\beta} \theta \cos^\alpha \phi \sin^\beta\phi (x + r \cos \phi \sin \theta)^a (y + r \sin \phi \sin \theta)^b 
\end{align}
Rewrite now
\begin{align}
 (x + r \cos \phi \sin \theta)^a (y + r \sin \phi \sin \theta)^b  = \sum_{k=0}^a \sum_{\ell=0}^b {a \choose k} {b \choose \ell} x^k y^\ell  r^{a-k+b-\ell} \sin^{b-\ell} \phi \cos^{a-k} \phi \sin^{a-k+b-\ell} \theta 
\end{align}

This proves the following
\begin{theorem}
With $\Delta \phi = \phi_\text{max} - \phi_\text{min}$, the spherical mean of monomial data reads
\begin{align}
 M^W[n_x^\alpha n_y^\beta x^ay^b](\vec x, r) &= \sum_{k=0}^a \sum_{\ell=0}^b {a \choose k} {b \choose \ell} x^k y^\ell  r^{a-k+b-\ell} 
 \frac{ \mu^W_{a+\alpha-k,b+\beta-\ell} }{2 \Delta \phi} \eta_{a+\alpha-k+b+\beta-\ell+1} \label{eq:sphericalmeanscalarfunction}
\end{align}
\end{theorem}

We rewrite it as
\begin{align}
 M^W[n_x^\alpha n_y^\beta x^ay^b](\vec x, r) &= \sum_{k=0}^a \sum_{\ell=0}^b F^{(a, b,\alpha,\beta)}_{k,\ell} x^k y^\ell  r^{a-k+b-\ell}  
\end{align}
with
\begin{align}
 F^{(a, b,\alpha,\beta)}_{k,\ell} := {a \choose k} {b \choose \ell} 
 \frac{ \mu^W_{a+\alpha-k,b+\beta-\ell} }{2 \Delta \phi} \eta_{a+\alpha-k+b+\beta-\ell+1}
\end{align}

In particular, if $x=y=0$, then the formula simplifies to
\begin{align}
 M^W[n_x^\alpha n_y^\beta x^ay^b](0, r) &= r^{a+b} \cdot \frac{ \mu^W_{a+\alpha,b+\beta} }{2 \Delta \phi} \eta_{a+\alpha+b+\beta+1}
\end{align}

\begin{corollary}
 If $\vec v = (x^a y^b, 0, 0)$, then the spherical mean $M^W[\vec n \cdot \vec v]$ is
 \begin{align}
  M^W[\vec n \cdot \vec v] &= M[\vec n x^a y^b]_x\\
  &= \frac{1}{2 \Delta \phi} \sum_{k=0}^a \sum_{\ell=0}^b {a \choose k} {b \choose \ell} x^k y^\ell  r^{a-k+b-\ell}  \mu_{a-k+1,b-\ell}(\phi_\text{max}, \phi_{min})  \eta_{a-k+b-\ell+2} \label{eq:sphericalmeanvectorxspr}
 \end{align}
 Analogously, if $\vec v = (0, x^a y^b, 0)$, then 
 \begin{align}
  M^W[\vec n \cdot \vec v] &= M[\vec n x^a y^b]_y\\
  &= \frac{1}{2 \Delta \phi} \sum_{k=0}^a \sum_{\ell=0}^b {a \choose k} {b \choose \ell} x^k y^\ell  r^{a-k+b-\ell} \mu_{a-k,b-\ell+1}(\phi_\text{max}, \phi_{min})   \eta_{a-k+b-\ell+2}
 \end{align}
\end{corollary}

\begin{corollary}
 If $\vec v = (x^a y^b, 0, 0)$, then the spherical mean $M[\vec n (\vec n \cdot \vec v)]$, restricted to two components, is
 \begin{align}
  M^W[\vec n (\vec n \cdot \vec v)] =\frac{1}{2 \Delta \phi} \sum_{k=0}^a \sum_{\ell=0}^b {a \choose k} {b \choose \ell} x^k y^\ell  r^{a-k+b-\ell} \vecc{ \mu_{a-k+2,b-\ell}(\phi_\text{max}, \phi_{min})  }{\mu_{a-k+1,b-\ell+1}(\phi_\text{max}, \phi_{min})  } \eta_{a-k+b-\ell+3} \label{eq:sphericalmeannnx}
 \end{align}
 If $\vec v = (0, x^a y^b, 0)$, then the spherical mean $M[\vec n \cdot \vec v]$, restricted to two components, is
 \begin{align}
  M^W[\vec n (\vec n \cdot \vec v)] =\frac{1}{2 \Delta \phi} \sum_{k=0}^a \sum_{\ell=0}^b {a \choose k} {b \choose \ell} x^k y^\ell  r^{a-k+b-\ell} \vecc{ \mu_{a-k+1,b-\ell+1}(\phi_\text{max}, \phi_{min})  }{\mu_{a-k,b-\ell+2}(\phi_\text{max}, \phi_{min})  } \eta_{a-k+b-\ell+3}
 \end{align}

\end{corollary}

\subsection{Initial data in $p$}

By linearity one can consider initial data in different variables separately.
First, consider initial data in $p$ only. The reconstructions in use are always polynomial; hence assume that $p_0(x, y) = x^a y^b$ in some wedge $W$. Then
\begin{align}
 p^W(t, \vec x) &= \del_r \left(r M[p_0](\vec x, r) \right) \Big|_{r=ct}\\
 &\overset{\eqref{eq:sphericalmeanscalarfunction}}{=} \sum_{k=0}^a \sum_{\ell=0}^b {a \choose k} {b \choose \ell} x^k y^\ell (a-k+b-\ell+1) (ct)^{a-k+b-\ell} 
 \frac{ \mu_{a-k,b-\ell}^W }{2 \Delta \phi} \eta_{a-k+b-\ell+1} \\
 &\overset{\vec x = 0}{=}(a+b+1) (ct)^{a+b}  \frac{ \mu_{ab}^W }{2 \Delta \phi} \eta_{a+b+1} \\
 \vec v^W(t, \vec x) &= - \frac{1}{r} \del_r \left( r^2 M[p_0  \vec n](\vec x, r) \right) \Big |_{r=ct}\\
 &\overset{\eqref{eq:sphericalmeanscalarfunction}}{=} -  \frac{1}{2 \Delta \phi} \sum_{k=0}^a \sum_{\ell=0}^b {a \choose k} {b \choose \ell} x^k y^\ell (a-k+b-\ell+2) (ct)^{a-k+b-\ell} \\\nonumber & \phantom{mmmmmmmmmmmmmmmmmmmmmmmm}\times\vecc{ \mu_{a-k+1,b-\ell}^W  }{\mu_{a-k,b-\ell+1}^W  } \eta_{a-k+b-\ell+2}\\
 &\overset{\vec x = 0}{=} -  \frac{1}{2 \Delta \phi}  (a+b+2) (ct)^{a+b} \vecc{ \mu_{a+1,b}^W  }{\mu_{a,b+1}^W  } \eta_{a+b+2}
\end{align}

\subsection{Initial data in $u$}

Consider now $\vec v_0 = (x^a y^b, 0, 0)$. Then, the solution reads
\begin{align}
 p^W(t, \vec x) &=  - \frac{1}{r} \del_r \left(r^2 M[\vec n \cdot \vec v_0](\vec x, r)\right ) \Big |_{r=ct}\\
 &\overset{\eqref{eq:sphericalmeanvectorxspr}}{=}  -  \frac{1}{2 \Delta \phi} \sum_{k=0}^a \sum_{\ell=0}^b {a \choose k} {b \choose \ell} x^k y^\ell (a-k+b-\ell+2) (ct)^{a-k+b-\ell}  \mu_{a-k+1,b-\ell}^W  \eta_{a-k+b-\ell+2} \\
 &\overset{\vec x = 0}{=} -  \frac{1}{2 \Delta \phi} (a+b+2) (ct)^{a+b}  \mu_{a+1,b}^W  \eta_{a+b+2}\\ 
 \vec v^W(t, \vec x) &= \frac23 \vec{ v}_0(\vec x)  + \del_r\left(rM[(\vec v_0 \cdot \vec n)  \vec n](\vec x, r) \right ) \nonumber \\&- M\left[  \vec v_0 - 3  (\vec v_0 \cdot \vec n)  \vec n   \right ](\vec x, r) - \int_0^{ct} \dd r \frac{1}{r} M \left[  \vec v_0 - 3  (\vec v_0 \cdot \vec n)  \vec n   \right ](\vec x, r)  \Big |_{r=ct}\\
 &\overset{\eqref{eq:sphericalmeanscalarfunction},\eqref{eq:sphericalmeannnx}}{=} \frac23 \vecc{x^a y^b}{0}  \\
 \nonumber &+ \frac{1}{2 \Delta \phi} \sum_{k=0}^a \sum_{\ell=0}^b {a \choose k} {b \choose \ell} x^k y^\ell (a-k+b-\ell+1) (ct)^{a-k+b-\ell} \vecc{ \mu_{a-k+2,b-\ell}^W  }{\mu_{a-k+1,b-\ell+1}^W  } \eta_{a-k+b-\ell+3}  \\ 
 \nonumber &-  \frac{1}{2 \Delta \phi}  \sum_{k=0}^a \sum_{\ell=0}^b {a \choose k} {b \choose \ell} x^k y^\ell  \vecc{   
  \mu_{a-k,b-\ell}^W \eta_{a-k+b-\ell+1} - 3 \mu_{a-k+2,b-\ell}^W\eta_{a-k+b-\ell+3} }{-3 \mu_{a-k+1,b-\ell+1}^W \eta_{a-k+b-\ell+3}} \\ \nonumber &\phantom{mmmmmmmmmmmmmmm} \times \left( (ct)^{a-k+b-\ell} +  \int_0^{ct} \dd r \,r^{a-k+b-\ell-1}  \right ) 
\end{align}

The integral $\int_0^{ct} \dd r \,r^{a-k+b-\ell-1}$ is
\begin{align}
 \frac{(ct)^{a-k+b-\ell}}{a-k+b-\ell} \qquad \text{if }a-k+b-\ell \neq 0
\end{align}
The condition is violated when $k=a$ and $\ell = b$. One can show that for continuous initial data, the term $\int_0^{ct} \dd r \frac{1}{r} M \left[  \vec v_0 - 3  (\vec v_0 \cdot \vec n)  \vec n   \right ](\vec x, r)$ is finite and in this case one can ignore the contribution from the case $a-k+b-\ell = 0$. In general, in particular for discontinuous data, a (logarithmic) singularity is part of the solution (see \cite{barsukow17}). Thus,
\begin{align}
 \vec v^W(t, \vec x) &= \frac23 \vecc{x^a y^b}{0}  \\
 \nonumber &+ \frac{1}{2 \Delta \phi} \sum_{k=0}^a \sum_{\ell=0}^b {a \choose k} {b \choose \ell} x^k y^\ell (ct)^{a-k+b-\ell} \\\nonumber & \phantom{mmm}\times\left [ (a-k+b-\ell+1)  \vecc{ \mu_{a-k+2,b-\ell}^W  }{\mu_{a-k+1,b-\ell+1}^W  } \eta_{a-k+b-\ell+3} \right . \\\nonumber&\phantom{mmmmm} \left. - \vecc{   
  \mu_{a-k,b-\ell}^W \eta_{a-k+b-\ell+1} - 3 \mu_{a-k+2,b-\ell}^W\eta_{a-k+b-\ell+3} }{-3 \mu_{a-k+1,b-\ell+1}^W \eta_{a-k+b-\ell+3}}  \left( 1 +  \frac{1}{a-k+b-\ell}  \right ) \right ]
\end{align}
where it is understood that the fraction in the last term is only to be inculded if it is defined.

If $\vec x = 0$, the expression simplifies:
\begin{align}
 \vec v^W(t, 0) &= \frac23 \vecc{x^a y^b}{0} + \frac{1}{2 \Delta \phi}  (ct)^{a+b} \\\nonumber & \phantom{mmm}\times\left [ (a+b+1)  \vecc{ \mu_{a+2,b}^W  }{\mu_{a+1,b+1}^W  } \eta_{a+b+3}  - \vecc{   
  \mu_{a,b}^W \eta_{a+b+1} - 3 \mu_{a+2,b}^W\eta_{a+b+3} }{-3 \mu_{a+1,b+1}^W \eta_{a+b+3}}  \left( 1 +  \frac{1}{a+b}  \right ) \right ]
\end{align}
where, again, the last fraction is only included if it is defined.

\section{Pseudocode implementation of the evolution operator} \label{app:evoop}

\begin{lstlisting}[escapeinside={/*@}{@*/},caption={Pseudocode implementation of the evolution operator.},captionpos=b]
FunctionPolynomial[] timeEvolution(FunctionPolynomial[] initial, double[] position, double /*@$\phi_\text{min}$@*/, double /*@$\phi_\text{max}$@*/){
/*@ This time evolution yields the result as a polynomial in $r$ that needs to be evaluated at $ct$@*/

    FunctionPolynomial[] res = new FunctionPolynomial[initial.length];
    double fraction = Math.abs(/*@$\phi_\text{max}$@*/ - /*@$\phi_\text{min}$@*/) / 2 / Math.PI;
			
    FunctionPolynomial R = new FunctionPolynomial(1, 0); 
    R.setCoefficient(1, 0, 1.0); // /*@polynomial $r$@*/
			
    /*@$M^W[p_0]$@*/ = sphericalMean(/*@$\phi_\text{min}$@*/, /*@$\phi_\text{max}$@*/, initial[PRESS], position);
			
    /*@$M^W[n_x p_0]$@*/ = sphericalMean(/*@$\phi_\text{min}$@*/, /*@$\phi_\text{max}$@*/, initial[PRESS], position, 1, 0);
    /*@$M^W[n_y p_0]$@*/ = sphericalMean(/*@$\phi_\text{min}$@*/, /*@$\phi_\text{max}$@*/, initial[PRESS], position, 0, 1);
			
    /*@$M^W[u_0]$@*/ = sphericalMean(/*@$\phi_\text{min}$@*/, /*@$\phi_\text{max}$@*/, initial[XVEL], position);
    /*@$M^W[v_0]$@*/ = sphericalMean(/*@$\phi_\text{min}$@*/, /*@$\phi_\text{max}$@*/, initial[YVEL], position);
			
    /*@$M^W[\vec n \cdot \vec v]$@*/ = sphericalMean(/*@$\phi_\text{min}$@*/, /*@$\phi_\text{max}$@*/, initial[XVEL], position, 1,0).add(
                     sphericalMean(/*@$\phi_\text{min}$@*/, /*@$\phi_\text{max}$@*/, initial[YVEL], position, 0,1) );

    /*@$M^W[n_x (\vec n \cdot \vec v)]$@*/ = sphericalMean(/*@$\phi_\text{min}$@*/, /*@$\phi_\text{max}$@*/, initial[XVEL], position, 2,0).add(
                     sphericalMean(/*@$\phi_\text{min}$@*/, /*@$\phi_\text{max}$@*/, initial[YVEL], position, 1,1) );
    /*@$M^W[n_y (\vec n \cdot \vec v)]$@*/ = sphericalMean(/*@$\phi_\text{min}$@*/, /*@$\phi_\text{max}$@*/, initial[XVEL], position, 1,1).add(
                     sphericalMean(/*@$\phi_\text{min}$@*/, /*@$\phi_\text{max}$@*/, initial[YVEL], position, 0,2) );
			
			
    PPterm = /*@$M^W[p_0]$@*/.mult(R).diffX();  
    PUVterm = (/*@$M^W[\vec n \cdot \vec v]$@*/.mult(R).mult(R).diffX()).divideByX();
    UPterm = ((/*@$M^W[n_x p_0]$@*/.mult(R).mult(R)).diffX()).divideByX();
    UUVterm = (((/*@$M^W[n_x (\vec n \cdot \vec v)]$@*/.mult(R).mult(R).mult(R).diffX()
                 ).divideByX().sub(/*@$M^W[u_0]$@*/.mult(R))).diffX()).divideByX().integrateX();
    VPterm = ((/*@$M^W[n_y p_0]$@*/.mult(R).mult(R)).diffX()).divideByX();
    VUVterm = (((/*@$M^W[n_y (\vec n \cdot \vec v)]$@*/.mult(R).mult(R).mult(R).diffX()
                ).divideByX().sub(/*@$M^W[v_0]$@*/.mult(R))).diffX()).divideByX().integrateX();
			
    res[PRESS] = PPterm.sub(PUVterm).mult(fraction);
    res[XVEL] = FunctionPolynomial.constantPolynomial(initial[XVEL].eval(position[0], position[1])
                                          ).sub(UPterm).add(UUVterm).mult(fraction);
    res[YVEL] = FunctionPolynomial.constantPolynomial(initial[YVEL].eval(position[0], position[1])
                                          ).sub(VPterm).add(VUVterm).mult(fraction);

    return res;
}


FunctionPolynomial sphericalMean(double /*@$\phi_\text{min}$@*/, double /*@$\phi_\text{max}$@*/, FunctionPolynomial /*@$f$@*/, double[] position, int /*@$\alpha$@*/, int /*@$\beta$@*/){
    /*@ Computes $M^W[n_x^\alpha y^\beta f]$ for polynomial $f$ at @*/position;
    /*@ Radius $r$ not specified, result of computation is polynomial in $r$ @*/
    FunctionPolynomial res = new FunctionPolynomial(/*@$f$@*/.degreeX() + /*@$f$@*/.degreeY()); // /*@degree of result@*/
		
    for (int i = 0; i <= /*@$f$@*/.degreeX(); i++){
        for (int j = 0; j <= /*@$f$@*/.degreeY(); j++){
            res = res.add(sphericalMeanMonomial(/*@$\phi_\text{min}$@*/, /*@$\phi_\text{max}$@*/, i, j, position, /*@$\alpha$@*/, /*@$\beta$@*/).mult(/*@$f$@*/.getCoefficient(i, j)));
        }
    }
    return res;
}
	
	
FunctionPolynomial sphericalMeanMonomial(double /*@$\phi_\text{min}$@*/, double /*@$\phi_\text{max}$@*/, int /*@a@*/, int /*@b@*/, double[] position, int /*@$\alpha$@*/, int /*@$\beta$@*/){
    /*@ Computes $M^W[n_x^\alpha y^\beta x^a y^b]$ at @*/position;
    /*@ Radius $r$ not specified, result of computation is polynomial in $r$ @*/
	
    FunctionPolynomial sphericalMean = new FunctionPolynomial(/*@a@*/ + /*@b@*/); // /*@degree of result@*/
		
    int phiIntervalIndex = /*@identify the quadrant required (can be passed directly as parameter)@*/
		
    for (int /*@$k$@*/ = 0; /*@$k$@*/ <= /*@$a$@*/; /*@$k$@*/++){
        for (int /*@$\ell$@*/ = 0; /*@$\ell$@*/ <=  /*@$b$@*/; /*@$\ell$@*/++){
            sphericalMean.addToCoefficient(/*@$a+b-\ell-k$@*/, 0, 
						factorInsideEvolution[/*@$a$@*/][/*@$b$@*/][/*@$k$@*/][/*@$\ell$@*/][/*@$\alpha$@*/][/*@$\beta$@*/][phiIntervalIndex]*
						power(position[0], /*@$k$@*/) * 
						power(position[1], /*@$\ell$@*/));
        }
    }
    return sphericalMean;
}


precomputeFactors(double /*@$\phi_\text{min}$@*/, double /*@$\phi_\text{max}$@*/, int phiIntervalIndex){
    for (int /*@$a$@*/ = 0; /*@$a$@*/ <= /*@$a_\text{max}$@*/; /*@$a$@*/++){
        for (int /*@$b$@*/ = 0; /*@$b$@*/ <= /*@$b_\text{max}$@*/; /*@$b$@*/++){
            for (int /*@$k$@*/ = 0; /*@$k$@*/ <= /*@$a$@*/; /*@$k$@*/++){
                for (int /*@$\ell$@*/ = 0; /*@$\ell$@*/ <= /*@$b$@*/; /*@$\ell$@*/++){
                    for (int /*@$\alpha$@*/ = 0; /*@$\alpha$@*/ <= /*@$\alpha_\text{max}$@*/; /*@$\alpha$@*/++){
                        for (int /*@$\beta$@*/ = 0; /*@$\beta$@*/ <= /*@$\beta_\text{max}$@*/; /*@$\beta$@*/++){
                            factorInsideEvolution[/*@$a$@*/][/*@$b$@*/][/*@$k$@*/][/*@$\ell$@*/][/*@$\alpha$@*/][/*@$\beta$@*/][phiIntervalIndex] = 
				    /*@${a \choose k} {b \choose \ell}$@*/ * sphericalMeanIntegrationFactor(/*@$\phi_\text{min}$@*/, /*@$\phi_\text{max}$@*/, /*@$a-k+\alpha$@*/, /*@$b-\ell+\beta$@*/);
			     /*@ This is $F^{(a, b,\alpha,\beta)}_{k,\ell}$. @*/
                        }
                    }
                }
            }
        }
    }
}

double sphericalMeanIntegrationFactor(double /*@$\phi_\text{min}$@*/, double /*@$\phi_\text{max}$@*/, int i, int j){
    /*@ This function computes $\eta_{i+j+1} \mu^W_{ij} / 2 / \Delta \phi$.@*/
    double secondFactor = 0.0;
    double factor = 0.0;

    double firstFactor = /*@$(i+j)!!/(i+j+1)!!$@*/;
    if ((i+j+1) /*@$\mod 2$@*/ == 0){
        firstFactor *= /*@$\pi$@*/;
    } else {
        firstFactor *= 2.0;
    }
    /*@ This is $\eta_{i+j+1}$.@*/
    

    for (int /*@$m$@*/ = 0; /*@$m$@*/ <= i; /*@$m$@*/++){
        for(int /*@$n$@*/ = 0; /*@$n$@*/ <= j; /*@$n$@*/++){
            factor = /*@${i \choose m} {j \choose n}$@*/;
            if (/*@$i+j-2m-2n$@*/ == 0){
                if(/*@$j\mod 2$@*/==0){
                    factor *= /*@$(-1)^{j/2+n} (\phi_\text{max}-\phi_\text{min})$@*/;
                } else {
                    factor = 0;
                }
            } else {
                factor *= /*@$(-1)^{n+j} / (i+j-2m-2n)$@*/;
                if(/*@$j\mod 2$@*/==0){
                    factor *= /*@$(-1)^{j/2}(\sin(\phi_\text{max}(i+j-2m-2n))-\sin(\phi_\text{min}(i+j-2m-2n)))$@*/;
                } else {
                    factor *= /*@$(-1)^{(j-1)/2} (\cos(\phi_\text{max}(i+j-2m-2n))-\cos(\phi_\text{min}(i+j-2m-2n)))$@*/;
                }
            }
            secondFactor += factor;
        }
    }
    secondFactor /= /*@$2^{i+j}$@*/; 

    return firstFactor*secondFactor / 2 / (/*@$\phi_\text{max}$@*/ - /*@$\phi_\text{min}$@*/);
}

\end{lstlisting}

The expressions are simplified tremendously through an implementation of polynomial functions; such an implementation is given in the Appendix \ref{app:details} as a reference, where the definitions of functions used in the above listings can be found.
It is possible to precompute not just some coefficients, but everything that is known a priori. This is best achieved algorithmically by using the previous implementation to compute necessary coefficients:

\begin{lstlisting}[escapeinside={/*@}{@*/},caption={Pseudocode implementation of the optimized evolution operator.},captionpos=b]
precomputeEvolutionCoefficients(){
    double[][] allPositions = /*@where one will want to evaluate the evolution operator, e.g. the 8 points at the boundary of cell@*/;
/*@  @*/
    /*@ Positions are assumed to be in $\left[-\frac{\Delta x}{2}, \frac{\Delta x}{2}\right] \times  \left[-\frac{\Delta y}{2}, \frac{\Delta y}{2}\right]$@*/
    
    evolutionCoefficients = new FunctionPolynomial[3][3][3][allRelativePositions.size()][4][3];
    /*@ The order is: variable of solution, $a$, $b$, position, quadrant, variable of initial data @*/
    /*@ A quadrant is the angle between@*/ /*@$\pi$@*//2*angleIndex /*@and@*/ /*@$\pi$@*//2*(angleIndex+1);
    
    FunctionPolynomial[] initial;
    int positionIndex;
		
    for (int q = 0; q <= 2; q++){ // variables u,v,p
        for (int /*@ a @*/ = 0; /*@ a @*/ <= 2; /*@ a @*/++){
            for (int /*@ b @*/ = 0; /*@ b @*/ <= 2; /*@ b @*/++){
                initial = new FunctionPolynomial[3];
                for (int q2 = 0; q2 <= 2; q2++){ initial[q2] = FunctionPolynomial.ZERO(); }
                initial[q] = new FunctionPolynomial(2, 2);
                initial[q].setCoefficient(/*@a@*/, /*@b@*/, 1.0);
                /*@ The q-th variable is now set to $x^a y^b$@*/
					
                positionIndex = 0;
                for (double[] position : allPositions) {
                    for (int angleIndex = 0; angleIndex <= 3; angleIndex++){
                        for (int q2 = 0; q2 <= 2; q2++){  // variables u,v,p
                            evolutionCoefficients[q][/*@a@*/][/*@b@*/][positionIndex][angleIndex][q2] = timeEvolution(initial, position, PI/2*angleIndex, PI/2*(angleIndex+1))[q2];
                        }
                    }
                    positionIndex++;
                }
            }
        }
    }
}
	
double[] timeEvolutionUsingPrecomputed(FunctionPolynomial[] initial, double[] position, double /*@$t$@*/, double /*@$\phi_\text{min}$@*/, double /*@$\phi_\text{max}$@*/){
	/*@ positions is assumed to be in $\left[-\frac{\Delta x}{2}, \frac{\Delta x}{2}\right] \times  \left[-\frac{\Delta y}{2}, \frac{\Delta y}{2}\right]$@*/
    double[] res = new double[initial.length];
    int positionIndex = /*@identify the index corresponding to@*/ position /*@(can be passed directly as parameter)@*/;
    int[] angleIndices = /*@identify the quadrants required (one or several; can be passed directly as parameter)@*/

    for (int q = 0; q <= 2; q++){ // variables u,v,p
        for (int /*@$a$@*/ = 0; /*@$a$@*/ <= 2; /*@$a$@*/++){
            for (int /*@$b$@*/ = 0; /*@$b$@*/ <= 2; /*@$b$@*/++){
                for (int angleIndex = 0; angleIndex < angleIndices.length; angleIndex++){
                    for (int q2 = 0; q2 <= 2; q2++){ // variables u,v,p
                        res[q2] += evolutionCoefficients[q][/*@$a$@*/][/*@$b$@*/][positionIndex][angleIndices[angleIndex]][q2].eval(/*@$ct$@*/) * initial[q].getCoefficient(/*@$a$@*/, /*@$b$@*/);
                    }
                }
            }
        }
    }
    return res;
}
\end{lstlisting}

\begin{lstlisting}[escapeinside={/*@}{@*/},caption={Example of the usage of the evolution operator.},captionpos=b]
if (edge.isVertical()){
    double[] leftValue  = timeEvolutionUsingPrecomputed(edge.leftCell().reconstruction(), new double[]{0.5*dx, 0.0}, /*@$\Delta t$@*/, /*@$\pi$@*//2, 3.0*/*@$\pi$@*//2);
    double[] rightValue = timeEvolutionUsingPrecomputed(edge.rightCell().reconstruction(), new double[]{-0.5*dx, 0.0}, /*@$\Delta t$@*/, -/*@$\pi$@*//2, /*@$\pi$@*//2);
    return leftValue + rightValue;
} else if (edge.isHorizontal(){
    double[] bottomValue = timeEvolutionUsingPrecomputed(edge.bottomCell().reconstruction(), new double[]{0.0, 0.5*dy}, /*@$\Delta t$@*/, -/*@$\pi$@*/, 0);
    double[] topValue    = timeEvolutionUsingPrecomputed(edge.topCell().reconstruction(), new double[]{0.0, -0.5*dy}, /*@$\Delta t$@*/, 0, /*@$\pi$@*/);
    return bottomValue + topValue;
}	
\end{lstlisting}

\section{Implementation details} \label{app:details}

\begin{lstlisting}[escapeinside={/*@}{@*/},caption={An example of an implementation of polynomial functions.},captionpos=b]
public class FunctionPolynomial extends Function{
	protected double[][] coefficients;

	public FunctionPolynomial(int degx, int degy) {
		coefficients = new double[degx+1][degy+1];
	}

	public void setCoefficient(int i, int j, double value){
		if ((i >= coefficients.length) || (j >= coefficients[0].length)){ resize(i, j); }
		coefficients[i][j] = value;
	}
	
	public FunctionPolynomial integrateX(){
		FunctionPolynomial res = new FunctionPolynomial(degreeX()+1, 0);
		res.setCoefficient(0, 0, 0.0);
		for (int j = 0; j < coefficients[0].length; j++){
			for (int i = 0; i < coefficients.length; i++){
				res.setCoefficient(i+1, j, coefficients[i][j]/(i+1)); 
			}
		}
		return res;
	}
	
	public FunctionPolynomial add(FunctionPolynomial poly){
		int[] otherDegree = poly.degree(); 
		int[] myDegree = this.degree();
		FunctionPolynomial res = new FunctionPolynomial(Math.max(otherDegree[0], myDegree[0]), Math.max(otherDegree[1], myDegree[1]));
		int[] deg = res.degree();
		
		for (int j = 0; j <= deg[1]; j++){
			for (int i = 0; i <= deg[0]; i++){
				res.setCoefficient(i, j, getCoefficient(i, j) + poly.getCoefficient(i,j)); 
			}
		}
		return res;
	}
	
	public static FunctionPolynomial ZERO(){ return constantPolynomial(0); }

	public static FunctionPolynomial constantPolynomial(double val){
		FunctionPolynomial res = new FunctionPolynomial(0, 0);
		res.setCoefficient(0, 0, val);
		return res;
	}

	[...]	
}

\end{lstlisting}

\end{document}